\documentclass[12pt]{article}

\usepackage{theorem}
\usepackage{amssymb}
\usepackage{amsmath}
\usepackage{amsfonts}
\usepackage{times}
\usepackage{graphicx}

\def \RR {\mathbb R}

\def \EE {\mathbb E}

\def \PP {\mathbb P}
\def \eps {\varepsilon}
\def \vphi {\varphi}

\def \cI {\mathcal I}

\def \cD {\mathcal D}
\def \cB {\mathcal B}

\def \cA {\mathcal A}

\newtheorem{theorem}{Theorem}
\newtheorem{lemma}{Lemma}

\newtheorem{corollary}{Corollary}
 {\theorembodyfont{\rmfamily}}

\def\myffrac#1#2 in #3{\raise 2.6pt\hbox{$#3 #1$}\mkern-1.5mu\raise 0.8pt\hbox{$
#3/$}\mkern-1.1mu\lower 1.5pt\hbox{$#3 #2$}}

\begin{document}

\title{A Berry-Esseen type inequality for convex bodies with an
unconditional basis}
\author{Bo'az Klartag\thanks{The author is a Clay Research Fellow,
and is also supported by NSF grant $\#DMS-0456590$. } \\ \\
{\normalsize Department of Mathematics} \\
{\normalsize Princeton University} \\
{\normalsize Princeton, NJ 08544, USA} \\
({\normalsize \it e-mail:} {\tt \normalsize bklartag@princeton.edu}) }
\date{}
\maketitle

\abstract{Suppose $X = (X_1,\ldots,X_n)$ is a random vector,
 distributed uniformly in a convex body $K \subset \RR^n$. We
assume the normalization $\EE X_i^2 = 1$ for $i=1,\ldots,n$.
The body $K$ is further required to be invariant under coordinate
reflections, that is, we assume that $(\pm X_1,\ldots,\pm X_n)$
has the same distribution as $(X_1,\ldots,X_n)$ for any choice of
signs. Then, we show that
$$
 \EE \left( \, |X| - \sqrt{n} \, \right)^2 \leq C^2,
$$
where $C \leq 4$ is a positive universal constant, and $| \cdot |$ is the
standard Euclidean norm in $\RR^n$. The estimate is tight, up to
the value of the constant. It leads to a Berry-Esseen type bound
in the central limit theorem for unconditional convex bodies. }


\section{Introduction}

Let $X_1,\ldots,X_n$ be random variables. We assume that the random
vector $X = (X_1,\ldots,X_n)$
is distributed according to a density $f: \RR^n \rightarrow [0, \infty)$,
and that the following hold:
\begin{enumerate}
\item[(A)] The joint density $f$ is log-concave. That is, the function $f$
has the form $f = e^{-H}$ with $H: \RR^n \rightarrow (-\infty, \infty]$
being a convex function.
\item[(B)] The joint density $f$ is ``unconditional''. That is, for any
point $(x_1,\ldots,x_n) \in \RR^n$ and a sign vector
$(\delta_1,\ldots,\delta_n) \in \{ \pm 1 \}^n$,
$$ f(x_1,\ldots,x_n) = f(\delta_1 x_1,\ldots, \delta_n x_n). $$
Equivalently, the random vector $(X_1,\ldots,X_n)$ has the same
distribution as $(\pm X_1,\ldots,\pm X_n)$ for any choice of
signs.
\item[(C)] The isotropic normalization $\EE X_i^2 = 1$ holds for $i=1,\ldots,n$.
\end{enumerate}
A particular case is when $X$ is distributed uniformly in a convex
set $K \subset \RR^n$, which is normalized so that $\EE X_i^2 = 1$
for all $i$, and is also ``unconditional'', i.e., for any $x =
(x_1,\ldots,x_n) \in \RR^n$ and for any choice of signs,
$$ (x_1,\ldots,x_n) \in K \ \ \ \ \ \Rightarrow \ \ \ \ \ \ (\pm
x_1,\ldots,\pm x_n) \in K. $$
We prove the following
Berry-Esseen type theorem:

\begin{theorem} Under assumptions (A), (B) and (C),
\begin{equation} \sup_{\alpha \leq \beta} \left| \PP \left( \alpha \leq
\frac{1}{\sqrt{n}} \sum_{i=1}^n X_i \leq \beta \right) -
\frac{1}{\sqrt{2 \pi}} \int_{\alpha}^{\beta} e^{-t^2 / 2} dt
\right|  \leq \frac{C}{n}, \label{eq_1124} \end{equation}
where $C
> 0$ is a universal constant. Moreover, for any
$\theta_1,\ldots,\theta_n \in \RR$ with $\sum_i \theta_i^2 = 1$,
\begin{equation}
 \sup_{\alpha \leq  \beta} \left| \PP \left( \alpha \leq \sum_{i=1}^n \theta_i X_i
 \leq \beta \right)
- \frac{1}{\sqrt{2 \pi}} \int_{\alpha}^{\beta} e^{-t^2 / 2} dt
\right|  \leq C \sum_{i=1}^n \theta_i^4. \label{eq_528}
\end{equation}
  \label{thm_511}
\end{theorem}

The log-concavity requirement (A) is crucial. A simple example may
be described as follows: Denote by $e_1,\ldots,e_n$ the standard
orthonormal basis in $\RR^n$. Let $T$ be a random variable,
 distributed uniformly in the set $\{1,\ldots,n\}$. Let $U$ be
 a random variable,
independent of $T$,
 distributed uniformly in the
interval $[-\sqrt{3 n}, \sqrt{3 n} ]$. Consider the random vector
$X = U e_T$. Then  $(\pm X_1,\ldots,\pm X_n)$ has the same
distribution as $(X_1,\ldots,X_n)$ for any choice of signs, and
also $\EE X_i^2 = 1$ for all $i$. However, $ \sum_i X_i =  U$ is
distributed uniformly in an interval, and hence its distribution
is far from normal. This demonstrates that assumptions (B) and (C)
alone cannot guarantee gaussian approximation.

\medskip
The bound in (\ref{eq_1124}) is optimal, up
to the precise value of the constant, as shown by the example of
$X_1,\ldots,X_n$ being independent random variables, with each
$X_i$ distributed, say, uniformly in a symmetric interval (see, e.g.,
\cite[Vol. II, Section XVI.4]{feller}).
A central element in the proof of Theorem \ref{thm_511} is
the sharp estimate
\begin{equation}
Var \left( \frac{|X|^2}{n} \right) = \EE \left( \frac{|X|^2}{n} - 1 \right)^2 \leq \frac{C}{n}, \label{eq_429_}
\end{equation}
for a positive universal constant $C \leq 16$.
Inequality
(\ref{eq_429_}) implies that most of the mass of the random
vector $X$ is concentrated in a {\it thin spherical shell}
of radius $\sqrt{n}$, centered at the origin in $\RR^n$, whose
width has the order of magnitude of a universal constant.
The bound (\ref{eq_429_}) was established by
Wojtaszczyk \cite{wjt}
in the case of Orlicz balls following a result
of Anttila, Ball and Perissinaki \cite{abp} regarding $\ell_p$-balls.
We say that a random vector $X = (X_1,\ldots,X_n)$ in $\RR^n$
is isotropically-normalized if
$\EE X_i = 0$ and $\EE X_i X_j = \delta_{i,j}$
for all $i,j$, where $\delta_{i,j}$ is Kronecker's delta.
A conjecture going back to
 Anttila, Ball and Perissinaki (see \cite{abp,BK}) is that
the thin spherical shell inequality (\ref{eq_429_}) actually holds whenever
$X$ is an isotropically-normalized random vector in $\RR^n$
with a log-concave density. We were
able to verify this conjecture under the additional assumption
that the density of $X$ is unconditional.

\medskip Theorem \ref{thm_511} ought to be understood
in the context of the
{\it central limit theorem for convex bodies}.
The central limit theorem for convex bodies
is the following high-dimensional effect,
suggested in
the works of Brehm and Voigt \cite{BV} and
 Anttila, Ball and Perissinaki \cite{abp}, and
proven in \cite{clt, power_law}:
Whenever $X = (X_1,\ldots,X_n)$ is an isotropically-normalized random vector in $\RR^n$, for large $n$,
with a log-concave density, then
for ``most'' choices of coefficients $\theta_1,\ldots,\theta_n \in
\RR$, the random variable $\sum_i \theta_i X_i$ is approximately
gaussian. (In the context of Theorem \ref{thm_511}, note
 that if the vector of coefficients
$(\theta_1,\ldots,\theta_n)$ is distributed uniformly on the unit sphere in $\RR^n$,
then the right-hand side of (\ref{eq_528}) is at most $C / n$ with probability
greater than $1 - C \exp(-c \sqrt{n})$. Here $C, c > 0$ are universal constants.)
There is an intimate relation between the central limit
theorem for convex bodies and thin spherical shell estimates
like (\ref{eq_429_}). This connection is well-known,
beginning with the work of Sudakov \cite{S}.
The reader is referred to, e.g., \cite{clt}
for more background on the
central limit theorem for convex bodies
and to, e.g., \cite{abp, B, BK} for the relation
to thin  shell estimates.

\medskip
Previous
techniques for obtaining
thin spherical shell estimates
under convexity
assumptions relied almost entirely on concentration of measure
ideas, either on the sphere (see \cite{fgp, clt}), or on the
orthogonal group (see \cite{power_law}). The quantitative
estimates that these techniques have yielded so far are
sub-optimal. Inequality
(\ref{eq_429_}) was previously known to hold
 with the bound $C / n^{\kappa}$ in place of $C / n$, where the exponent
$\kappa$ is slightly smaller than $1/5$, see
\cite{clt,power_law}. The latter result is applicable for all isotropically-normalized random
vectors with a log-concave density.

\medskip In this article we suggest a different approach.
 Rather
than employing concentration of measure inequalities, our proof of
the optimal inequality (\ref{eq_429_}) is based on
analysis of the Neumann Laplacian on convex domains, the so-called
$L^2$-method in convexity, going back
to H\"ormander \cite{hormander} and to Helffer and Sj\"ostrand
\cite{HS}. The argument is further simplified by
using the theory of optimal
transportation of measures. We
expect this technique to be useful also in the study of other
problems in convex geometry, such as central limit theorems for
convex bodies with various types of symmetries.
The argument leading to the thin shell estimate occupies Section \ref{sec_laplace},
Section \ref{sec_transport} and Section \ref{sec_uncond}.
In Section
\ref{sec_berry} we apply these estimates and complete the proof of Theorem \ref{thm_511}.

\medskip Readers who are interested only in the proof of  inequality
(\ref{eq_429_}) and Theorem \ref{thm_511} may skip Section
\ref{sec_digression}. This section
is devoted to several results,
that were obtained as by-products, regarding the first non-zero
eigenvalue and the corresponding eigenfunctions of the Neumann
Laplacian on $n$-dimensional convex bodies. In particular, we show
 that the eigenfunctions are all ``biased'' towards some
direction in space. This rules out, for instance, the possibility
of an even eigenfunction.

\medskip
As the reader has probably figured out by now, we denote
expectation by $\EE$ and probability by $\PP$. We write $Var$ for
variance, and $Vol_n(A)$ for the Lebesgue measure of a measurable
set $A \subset \RR^n$. The scalar product of $u,v \in \RR^n$ is
denoted by $u \cdot v$. The letters $c, C, C^{\prime}, \tilde{c}$
etc. stand for various positive universal constants, whose value
may change from one line to the next.

\medskip
\emph{Acknowledgement.} We would like to express our gratitude to
Sasha Sodin for his kind help with the analysis related
to the classical central limit theorem, to Tom Spencer for illuminating
explanations regarding the work of
Helffer and Sj\"ostrand, and to Dario Cordero-Erausquin, Leonid Friedlandler,
Robert McCann, Emanuel Milman, Vitali Milman and Elias Stein for valuable
discussions on related topics. Thanks also to the referee for useful
comments and suggestions.

\section{Convexity and the Neumann Laplacian}
\label{sec_laplace}

In this section we analyze some convexity related properties of the Neumann
Laplacian, most of which are standard.  
 A convex body in $\RR^n$ is a compact,
convex set with a non-empty interior.
Let $K \subset \RR^n$ be a convex body
with a $C^{\infty}$-smooth boundary, to be fixed throughout this
section. We say that a function $\vphi: K \rightarrow \RR$ belongs
to $C^{\infty}(K)$ if all of its derivatives of all
orders exist and are bounded in the interior of $K$. When $\vphi$
is a $C^{\infty}(K)$-smooth function, the boundary
values of $\vphi$ and its derivatives are well defined, and are
$C^{\infty}$-smooth on the boundary $\partial K$.
 For $u \in
C^{\infty}(K)$  define
$$ \| u \|_{H^{-1}(K)} = \sup \left \{ \int_{K} \vphi
u \,  ; \, \vphi \in C^{\infty}(K), \ \int_K |\nabla \vphi|^2
\leq 1 \right \}. $$
Note that necessarily $\| u \|_{H^{-1}(K)} = \infty$ when $\int_K u \neq 0$.
For a
function $f$ in $n$ variables
and for $i=1,\ldots,n$ we write $\partial^i f$ for the derivative
of $f$ with respect to the $i^{th}$ coordinate.
When $f: K \rightarrow \RR$ is a square-integrable function, set
$$ Var_K(f) = \int_K \left( f(x) - E \right)^2 dx $$
with $E = Vol_n(K)^{-1} \int_K f$.
The main result of this section reads as follows:

\begin{lemma} Let $K \subset \RR^n$ be a convex body
with a  $C^{\infty}$-smooth boundary. Let $f: K
\rightarrow \RR$ be a $C^{\infty}(K)$-smooth function.
Then,
\begin{equation}
 Var_K(f) \leq \sum_{i=1}^n \| \partial^i f \|_{H^{-1}(K)}^2.
 \label{eq_211}
\end{equation}  \label{main_laplace}
\end{lemma}

One may verify that the right-hand side of (\ref{eq_211}) does not
depend on the choice of orthogonal coordinates in $\RR^n$.
See \cite{vil} for an analog of Lemma
\ref{main_laplace} for non-convex domains.
Let
$\rho: K \rightarrow \RR$ be a convex function which is
$C^{\infty}$-smooth with bounded derivatives of all orders in a
neighborhood of $\partial K$, such that
$$ \rho(x) = 0, \ |\nabla \rho(x)| = 1 \ \ \ \ \ \text{for} \ x \in \partial K $$
and $\rho(x) \leq 0$ for $x \in K$. For instance, we may select
$\rho(x) = - d(x,
\partial K) =  - \inf_{y \in
\partial K} |x - y|$. Note that for any $x \in \partial K$, the
vector $\nabla \rho(x)$ is the outer unit normal
to $\partial K$ at  $x$.

\medskip
Denote by $\cD$ the space of all
$C^{\infty}(K)$-smooth functions $u: K \rightarrow \RR$
that satisfy the following Neumann boundary condition: $$ \nabla u(x) \cdot \nabla \rho(x) = 0 \ \ \ \ \ \text{for} \ x
\in \partial K. $$
The following lemma
is a standard Bochner-Weitzenb\"ock
type integration by parts formula,
going back at least to
Lichnerowicz \cite{Lich}, to H\"ormander  \cite{hormander}
and to Kadlec \cite{kadlec}.
We write $\nabla^2 u$ for the hessian matrix of the function $u$.

\begin{lemma} Let $u \in \cD$ and denote $f = -\triangle u$.
Then,
\begin{equation}
 \int_K f^2 =
 \int_K  \sum_{i=1}^n|\nabla \partial^i
u|^2 + \int_{\partial K} \nabla^2 \rho (\nabla u) \cdot \nabla u.
\label{hormander}
\end{equation}
\label{lem_159}
\end{lemma}

\emph{Proof:} The function $x \mapsto \nabla
u(x) \cdot \nabla \rho(x)$ vanishes on $\partial K$. Since $\nabla
u$ is tangential to $\partial K$, the derivative of the function
$x \mapsto \nabla
u(x) \cdot \nabla \rho(x)$ in the direction
of $\nabla u$ vanishes on $\partial K$. That is,
$$ \nabla u(x) \cdot \nabla \left( \nabla u(x) \cdot \nabla \rho(x) \right) = 0 \ \ \ \ \ \text{for} \ \  x \in \partial K.
$$
Equivalently,
\begin{equation}
 (\nabla^2 u) (\nabla \rho) \cdot \nabla u + (\nabla^2 \rho) (\nabla
u) \cdot \nabla u = 0 \ \ \ \ \ \ \ \text{on}  \ \ \partial K.
\label{eq_1207}
\end{equation}
By Stokes theorem,
\begin{equation}
\int_K f^2 = \int_K (\triangle u)^2 = -\int_K \nabla (\triangle u) \cdot
\nabla u + \int_{\partial K} (\triangle u  \nabla u) \cdot
\nabla \rho. \label{eq_545}
\end{equation}
The boundary term vanishes, since $\nabla u \cdot
\nabla \rho = 0$ on $\partial K$. We conclude from
(\ref{eq_545}) and from an additional application of Stokes theorem that
$$ \int_K f^2 = -\sum_{i=1}^n \int_K
 \partial^i u
\triangle(\partial^i u)
= \sum_{i=1}^n \int_K |\nabla \partial^i u|^2 -
 \int_{\partial K} \sum_{i=1}^n (\partial^i u \nabla
\partial^i u ) \cdot \nabla \rho.
$$
Note that the integrand in the integral over $\partial K$ is
exactly $\nabla^2 u (\nabla \rho) \cdot \nabla u$. Hence, from
(\ref{eq_1207}),
$$ \int_K f^2 = \sum_{i=1}^n \int_K |\nabla \partial^i u|^2 + \int_{\partial K}
\nabla^2 \rho(\nabla u) \cdot \nabla u, $$ and the lemma is
proven. \hfill $\square$

\medskip The convexity of $K$ will be used next.
Recall that $\rho$ is a convex function, and hence its
hessian $\nabla^2 \rho(x)$ is a positive semi-definite matrix for
any $x \in \partial K$. Therefore, Lemma \ref{lem_159} implies
that for any $u \in \cD$,
\begin{equation}
  \sum_{i=1}^n \int_K |\nabla \partial^i
u|^2
\leq \int_K f^2
\label{eq_600}
\end{equation}
where $f = \triangle u$.
Lemma \ref{main_laplace} will be proven by dualizing  inequality
(\ref{eq_600}), in a way which is very much related to the approach taken by
H\"ormander \cite{hormander} and by Helffer and
Sj\"ostrand \cite{HS}.

\medskip \emph{Proof of Lemma \ref{main_laplace}:}
We are given $f \in C^{\infty}(K)$ and
we would like to prove  (\ref{eq_211}). We may assume that $\int_K
f = 0$ (otherwise, subtract $\frac{1}{Vol_n(K)} \int_K f$ from the function
$f$).

\medskip
Since $f \in C^{\infty}(K)$ and $\int_K  f = 0$,
there exists $u \in \cD$ with
$$ -\triangle u = f. $$
The existence of such $u \in \cD$ is a consequence of the classical existence and regularity theory of the
Neumann problem for the Laplacian on domains with a
$C^{\infty}$-smooth boundary (see, e.g.,
Folland's book \cite[chapter 7]{folland}).
Stokes
theorem yields
$$
 \int_K f^2 = -\int_K f \triangle u = \int_K \nabla f \cdot \nabla u - \int_{\partial K} f \nabla u \cdot \nabla \rho
= \sum_{i=1}^n \int_K \partial^i f \partial^i u,
$$
where the boundary term vanishes since $u \in \cD$.
From the definition of the $H^{-1}(K)$-norm and the
Cauchy-Schwartz inequality,
\begin{eqnarray} \label{eq_659}
 \int_K f^2 = \sum_{i=1}^n \int_K \partial^i f \partial^i u
& \leq & \sum_{i=1}^n \| \partial^i f \|_{H^{-1}(K)} \sqrt{\int_K |\nabla \partial^i u|^2} \\ & \leq & \sqrt{ \sum_{i=1}^n \| \partial^i f \|_{H^{-1}(K)}^2 } \cdot
\sqrt{ \sum_{i=1}^n \int_K |\nabla \partial^i u|^2 }. \nonumber
\end{eqnarray}
Combine (\ref{eq_659}) and (\ref{eq_600}) to conclude that
 $$ \int_K f^2 \leq \sum_{i=1}^n \| \partial^i f \|_{H^{-1}(K)}^2. $$
\hfill $\square$

\section{Transportation of Measure}
\label{sec_transport}

Suppose $\mu_1$ and $\mu_2$ are finite Borel measures on $\RR^m$
and $\RR^n$ respectively, and $T: \RR^m \rightarrow \RR^n$ is a
measurable map. We say that $T$ pushes forward, or transports,
$\mu_1$ to $\mu_2$ if
$$ \mu_1( T^{-1} (A) ) = \mu_2( A ) $$
for all Borel sets $A \subseteq \RR^n$. In this case we write
$\mu_2 = T_{\#} \mu_1$, and we call $T$ the transportation map.
Note that $\int ( \vphi \circ T ) d\mu_1 = \int \vphi d (T_{\#}
\mu_1)$ for any bounded, measurable function $\vphi$.

\medskip For example, let $\gamma$ be a Borel measure on $\RR^n
\times \RR^n$. For $(x,y) \in \RR^n \times \RR^n$ we write
$P^1(x,y) = x$ and $P^2(x,y) = y$. We say that the measure $P^1_{\#} \gamma$
is the marginal of $\gamma$ on the first coordinate, and $P^2_{\#}
\gamma$ is the marginal of $\gamma$ on the second coordinate. A
measure $\gamma$ on $\RR^n \times \RR^n$ with $P^1_{\#} \gamma =
\mu_1$ and $P^2_{\#} \gamma = \mu_2$ is called a ``coupling'' of
$\mu_1$ and $\mu_2$.

\medskip Suppose $\mu_1$ and $\mu_2$ are two finite Borel
measures on $\RR^n$. If $T$ pushes forward $\mu_1$ to $\mu_2$, then the map
$$ x \mapsto (x, T x)  $$
transports the measure $\mu_1$ to a measure $\gamma$ on $\RR^n
\times \RR^n$ which is a coupling of $\mu_1$ and $\mu_2$. The
$L^2$-Wasserstein distance between $\mu_1, \mu_2$ is defined as
$$ W_2(\mu_1, \mu_2) = \inf_{\gamma} \left( \int_{\RR^n \times \RR^n} |x - y|^2 \, d \gamma(x,y) \right)^{1/2}, $$
where the infimum runs over all couplings $\gamma$ of $\mu_1$ and
$\mu_2$. If there is no coupling, then $W_2(\mu_1, \mu_2) =
\infty$.
Let $\mu$ be a finite, compactly-supported Borel measure on $\RR^n$.
For a $C^{\infty}$-smooth function $u: \RR^n \rightarrow \RR$, set
$$ \| u \|_{H^{-1}(\mu)} = \sup \left \{ \int_{\RR^n} u \vphi
 \, d\mu  \, ; \, \vphi \in C^{\infty}(\RR^n), \ \int_{\RR^n}
|\nabla \vphi|^2 \, d\mu \leq 1 \right \}. $$ This definition fits
with the one given in Section \ref{sec_laplace}; We have $\| u
\|_{H^{-1}(\lambda_K)} = \| u \|_{H^{-1}(K)}$ where $\lambda_K$
denotes the restriction of the Lebesgue measure to $K$.

\medskip The next theorem is an
 extension of a remark by Yann Brenier \cite{brenier}
that we learned from Robert McCann.
For the convenience of the reader, we provide in the appendix
a detailed exposition of the elegant proof from
Villani \cite[Section 7.6]{villani}.

\begin{theorem} Let $\mu$ be a finite, compactly-supported
Borel measure on $\RR^n$. Let $h: \RR^n \rightarrow \RR$ be
a bounded, measurable function with $$ \int h d \mu = 0. $$ For a
sufficiently small $\eps > 0$, let $\mu_{\eps}$ be the measure whose
density with respect to $\mu$ is the non-negative function $1 +
\eps h$. Then,
$$ \| h \|_{H^{-1}(\mu)} \leq \liminf_{\eps \rightarrow 0^+} \frac{W_2(\mu, \mu_{\eps})}{\eps}.
 $$
\label{thm_258}
\end{theorem}

See \cite{brenier} and \cite{villani} for the intuition behind Theorem \ref{thm_258}.
We write $e_1,\ldots,e_n$ for the standard orthonormal basis in $\RR^n$.
Let $K \subset \RR^n$ be a convex body. Fix a point $x \in K$ and
$i=1,\ldots,n$. Consider the line $x + \RR e_i$, that is, the line in the direction of
$e_i$ that passes through $x$. This line meets $K$ with a
closed segment (or a single point). The two endpoints of this
segment in $\RR^n$ will be denoted by $\cB_i^-(x)$ and
$\cB_i^+(x)$, where $\cB_i^-(x) \cdot e_i \leq \cB_i^+(x) \cdot e_i$. Thus,
$$ K \cap \left( x + \RR e_i \right) = [\cB_i^-(x), \cB_i^+(x)], $$
the line segment from $\cB_i^-(x)$ to $\cB_i^+(x)$. See Figure 1.

\medskip For $i=1,\ldots,n$ consider the projection
 $$ \pi_i(x_1,\ldots,x_n) =  (x_1,\ldots,x_{i-1}, x_{i+1},\ldots,x_n), $$
defined for $(x_1,\ldots,x_n) \in \RR^n$.
Then $\pi_i(K)$ is a convex body in $\RR^{n-1}$.
For $y \in \pi_i(K)$, we define $q_i^{-}(y) \in \RR$
to be the minimal $i^{th}$ coordinate among all points $x \in K$ with
$\pi_i(x) = y$. Similarly, we define $q_i^+(y)$
to be the maximal $i^{th}$ coordinate.

\begin{center}
\includegraphics[height=1.5in]{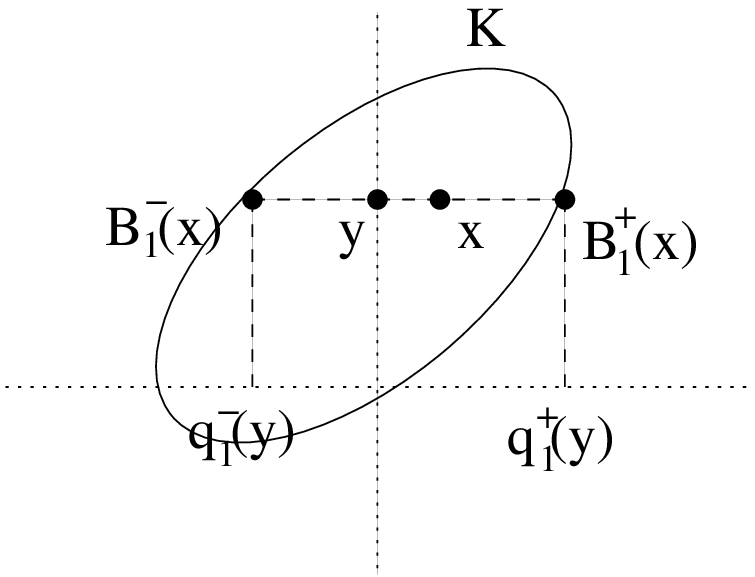}

{\small Figure 1}
\end{center}

\begin{lemma} Let $K \subset \RR^n$ be a convex body with a $C^{\infty}$-smooth
boundary. Fix $i=1,\ldots, n$. Let $\Psi: K \rightarrow \RR$ be a
$C^{\infty}(K)$-smooth function such that for any $x  \in K$,
\begin{equation}
 \Psi \left( \cB_i^-(x) \right) = \Psi \left( \cB_i^+ (x) \right).
 \label{eq_1231}
\end{equation}
For a sufficiently small $\eps > 0$ denote by $\mu_{\eps}$ the
measure whose density with respect to $\mu$ is $1 + \eps
\partial^i \Psi$. Then,
$$ \liminf_{\eps \rightarrow 0^+} \frac{W_2(\mu, \mu_{\eps})}{\eps} \leq
\sqrt{\int_K \left[ \Psi(x) - \Psi(\cB^+_i(x))
 \right]^2 dx}.
$$ \label{lem_1117}
\end{lemma}

\emph{Proof:} Without loss of generality, assume that $i = 1$.
For a sufficiently small $\eps > 0$, the function $1 + \eps \partial^1 \Psi$ is
positive on $K$, and hence $\mu_{\eps}$ is a non-negative measure.
Fix such a sufficiently small $\eps > 0$.

\medskip For $x = (t,x_2,\ldots,x_n) \in \RR^n$ we will use
the coordinates $x = (t, y)$ where $y = (x_2,\ldots,x_n) \in \RR^{n-1}$.
Fix $y \in \pi_1(K)$ and
denote $p = q^-_1(y)$ and $q =
q^+_1(y)$. According to our assumption
(\ref{eq_1231}),
$$  \int_p^q (1 + \eps \partial^1 \Psi(t,y) ) dt = (q - p) + \left. \eps
\Psi(t,y) \right|_{t=p}^q = q - p.
$$
Consequently, the densities $t \mapsto 1$ and  $t \mapsto
1 + \eps
\partial^1 \Psi(t,y)$ have an equal amount of mass on the interval
$[p,q]$. We consider the monotone transportation between these two
densities. That is, we define a map $T = T^{y}: [p,q]
\rightarrow [p,q]$ by requiring that for any $x_1 \in [p,q]$,
\begin{equation}
\int_{p}^{x_1} \left(
 1 + \eps \partial^1 \Psi(t,y) \right) dt =  \int_{p}^{T(x_1)} dt.
\label{eq_1237}
\end{equation}
The unique map $T: [p,q] \rightarrow [p,q]$ that satisfies
(\ref{eq_1237}) transports the measure whose density is $1 +
\eps
\partial^1 \Psi(t,y)$ on $[p,q]$ to the Lebesgue
measure on $[p,q]$. We deduce from (\ref{eq_1237}) that for $x_1
\in [p, q]$,
$$ T(x_1) = x_1 + \eps \left[ \Psi(x_1,y) - \Psi(p,y) \right]. $$
Therefore, \begin{eqnarray} \label{eq_105} \lefteqn{ \int_p^q | T(t) -
t |^2 \cdot
\left( 1 + \eps \partial^1 \Psi(t,y) \right) dt } \\
&  = & \eps^2 \int_p^q \left[ \Psi(t,y) -
\Psi(p,y) \right]^2 dt + \eps^3 R, \nonumber
\end{eqnarray}
with $|R|$ bounded by a constant depending only on $\Psi$ and $K$
(and in particular, independent of $\eps$ or $y$). We
now let $y \in \pi_1(K)$ vary, and we write
$$ S(x_1,y) = \left(  T^{y}(x_1), y \right)
\ \ \ \ \ \ \ \ \ \ \ \ \ \ \ \ \text{for} \ \ \  (x_1,y) \in K. $$ Note that $S$ is
well-defined (since $x_1$ belongs to the domain of
definition of $T^{y}$ when $(x_1,y) \in
K$), one-to-one, continuous, and maps $K$ onto $K$. Moreover, by Fubini, for
any continuous function $\vphi: K \rightarrow \RR$,
\begin{eqnarray*}
\lefteqn{ \int_K \vphi(S(x)) d \mu_{\eps}(x)  \, \, \, = \, \,
 \int_{\pi(K)} \left[
\int_{q_1^-(y)}^{q_1^+(y)}
\vphi(T^{y}(x_1), y) \cdot (1 + \eps \partial^1 \Psi) dx_1 \right] dy } \\
& =& \int_{\pi(K)} \left[
\int_{q_1^-(y)}^{q_1^+(y)} \vphi(x_1,
y) dx_1 \right] dy = \int_K \vphi(x) d
\mu(x). \phantom{aaaaaaaaaaaaaa}
\end{eqnarray*}
Therefore the map $S$ transports $\mu_{\eps}$ to $\mu$. According
to (\ref{eq_105}),
$$
 W_2(\mu, \mu_{\eps})^2 \leq \int_K |S(x) - x|^2 d
\mu_{\eps}(x) = \eps^2 \int_K \left[ \Psi(x) - \Psi \left(
\cB_1^-(x) \right) \right]^2 dx + \eps^3 R^{\prime},
$$
with $|R^{\prime}|$ smaller than a constant depending
only on $K$ and $\Psi$, and in particular independent of $\eps$.
To complete the proof, let $\eps$ tend to zero. \hfill $\square$

\section{A digression: Neumann eigenvalues and eigenfunctions}
\label{sec_digression}

This section presents some additional relations between convexity
and the Neumann Laplacian.
We retain the setup and notation of
Section \ref{sec_laplace}. We write $L^2(K)$ for the Hilbert space
that is the completion of $C^{\infty}(K)$ with respect to the norm
$$ \| u \|_{L^2(K)} = \sqrt{\int_K u^2}. $$ The operator
$-\triangle$, acting on the subspace $\cD \subset L^2(K)$, is a
symmetric, positive semi-definite operator. The classical theory implies that $-\triangle$ has
a complete system of orthonormal Neumann eigenfunctions
$\vphi_0,\vphi_1,\ldots \in \cD$ and Neumann eigenvalues
$0 \leq \lambda_0 \leq \lambda_1 \leq ...$ (see, e.g., \cite[Chapter
7]{folland}). The first eigenvalue is $\lambda_0 = 0$, with the
eigenfunction $\vphi_0$ being constant. It is well-known that
$\lambda_1 > 0$ when $K$ is convex (see, e.g, \cite{pw}. It is actually
enough to assume that $K$ is connected, see e.g., \cite[Theorem 1]{chavel_spectral}). We
refer to $\lambda_1$ as the first non-zero Neumann eigenvalue of $K$. It is well-known
that for any
$C^{\infty}(K)$-smooth function $u$ with $\int_K u = 0$,
\begin{equation}\lambda_1 \int_K u^2 \leq \int_K |\nabla u|^2.
\label{eq_131}
\end{equation}
Equality in (\ref{eq_131}) holds if and only if $u$ is an eigenfunction
corresponding to the eigenvalue $\lambda_1$.

\medskip
We say that the boundary of $K$ is uniformly strictly convex if
$\nabla^2 \rho(x)$ is a positive definite matrix for any $x \in
\partial K$. Equivalently, $\partial K$ is uniformly strictly
convex if the principal curvatures are all positive --  and not
merely non-negative -- everywhere on the boundary. Our next
corollary claims, loosely speaking, that  any non-trivial
eigenfunction corresponding to $\lambda_1$ cannot be ``spatially
isotropic'', but must have ``preference'' for a certain direction
in space.

\begin{corollary} Suppose $K \subset \RR^n$ is a convex body whose boundary
is $C^{\infty}$-smooth and uniformly strictly convex. Let $0 \not
\equiv \vphi \in \cD$ be an eigenfunction corresponding to the
first non-zero Neumann eigenvalue. Then,
\begin{equation} \int_K \nabla \vphi \neq 0.
\label{eq_334}
\end{equation}
Consequently, the multiplicity of the first non-zero Neumann
eigenvalue is at most $n$. \label{cor_1154}
\end{corollary}

\emph{Proof:} Assume the opposite. Then,
\begin{equation}
 \int_K \partial^i \vphi = 0 \ \ \ \ \text{for} \ i=1,\ldots,n.
 \label{eq_1113_}
\end{equation}  We write $\lambda_1$ for the first
non-zero eigenvalue, i.e., $\triangle \vphi = -\lambda_1 \vphi$. Since $\vphi \in \cD$,
inequality (\ref{eq_600}) gives
\begin{equation} \label{eq_1113}
\lambda_1^2 \int_K \vphi^2 = \int_K |\triangle \vphi|^2
 \geq  \sum_{i=1}^n \int_K |\nabla
\partial^i \vphi|^2.
\end{equation}
From (\ref{eq_1113_})
we know that $\int_K \partial^i \vphi = 0$ for all $i$. Thus (\ref{eq_1113}) and (\ref{eq_131}) yield
$$ \lambda_1^2 \int_K \vphi^2 \geq  \sum_{i=1}^n \int_K |\nabla
\partial^i \vphi|^2 \geq \lambda_1 \sum_{i=1}^n
\int_K (\partial^i \vphi)^2 = \lambda_1 \int_K |\nabla \vphi|^2 = \lambda_1^2
\int_K \vphi^2. $$
Therefore, there must be equality in all steps
and hence $\partial^1 \vphi,\ldots,\partial^n \vphi$ are
all Neumann eigenfunctions with eigenvalue $\lambda_1$.
We necessarily have equality also in
(\ref{eq_1113}). According to Lemma \ref{lem_159} this means that
$$
\int_{\partial K} \nabla^2 \rho(\nabla \vphi) \cdot \nabla \vphi =
0. $$ Since the integrand is non-negative and continuous, necessarily
\begin{equation} \label{eq_1130}
\nabla^2 \rho(\nabla \vphi) \cdot \nabla \vphi = 0 \ \ \ \ \ \
\text{on} \ \ \partial K.
\end{equation}

\medskip So far we have only used the convexity of $K$. The uniform
strict convexity of $\partial K$ means that $\nabla^2 \rho > 0$
on $\partial K$. Equation  (\ref{eq_1130})  has the
consequence that $\nabla \vphi = 0$ on $\partial K$, and therefore
\begin{equation} \vphi \equiv Const \ \ \ \ \ \ \text{on} \ \ \partial
K. \label{eq_1059} \end{equation}
This is well-known to be impossible for a Neumann eigenfunction
corresponding to the first non-zero eigenvalue. We sketch the standard
argument, see, e.g., \cite{chavel_spectral} for more information. Denote
$$ N=\{ x \in K ; \vphi(x) > 0 \}.$$ The set $N$ is non-empty since $\int_K
\vphi = 0$. Moreover, $\vphi$ vanishes on $\partial N$ because of (\ref{eq_1059}).
Since $\triangle \vphi = -\lambda_1 \vphi$ in $N$, then $\vphi$ is a Dirichlet
eigenfunction of the domain $N$ corresponding to the Dirichlet eigenvalue $\lambda_1$.
For a domain $\Omega \subset \RR^n$, denote by $\lambda_0^D(\Omega)$
the minimal eigenvalue of $-\triangle$ with Dirichlet boundary conditions
on $\Omega$. Then $\lambda_0^D(N) \leq \lambda_1$, as is witnessed by $\vphi$. Furthermore,
$\lambda_0^D(N) \geq \lambda_0^D(K)$ by domain monotonicity (see, e.g, \cite{chavel_spectral}),
hence $\lambda_0^D(K) \leq \lambda_1$.
However, we have the strict inequality $\lambda_0^D(K) > \lambda_1$
 (see, e.g., \cite{lw} for a much more accurate result).
We thus arrive at a contradiction. Consequently our assumption that $\int_K \nabla \vphi = 0$
was absurd. The proof of (\ref{eq_334}) is complete.

\medskip The linear map $\vphi \mapsto \int_K \nabla \vphi$
from the eigenspace of $\lambda_1$ to $\RR^n$ is therefore injective,
so the multiplicity of the eigenvalue cannot exceed $n$. \hfill $\square$

\medskip \emph{Remark.} Leonid Friedlandler explained to us how to eliminate the uniform strict convexity requirement from Corollary \ref{cor_1154}. 
His idea is to observe that since 
$\partial^1 \vphi,\ldots, \partial^n \vphi$ are all eigenfunctions,
then the restriction of $\vphi$ to the boundary $\partial K$ is actually an eigenfunction of the Laplacian associated with the Riemannian manifold $\partial K$.
However, (\ref{eq_1130}) entails that $\vphi$ is constant in some open 
set in $\partial K$, which is known to be impossible for an eigenfunction. We omit the details.

\medskip For $i=1,\ldots,n$ and $x = (x_1,\ldots,x_n) \in \RR^n$
write
$$ \sigma_i(x) = (x_1,\ldots,x_{i-1}, -x_i, x_{i+1}, \ldots, x_n),
$$
i.e., we flip the sign of the $i^{th}$ coordinate. For a function
$f$, we write $\sigma_i(f)(x) = f(\sigma_i(x))$. Our
next corollary exploits the well-known relationship between the
eigenfunctions and symmetry. Similar arguments appear, e.g., in
\cite{BB}.

\begin{corollary} Suppose $K \subset \RR^n$ is a convex body with
a $C^{\infty}$-smooth boundary.  Denote by $E_{\lambda_1} \subset \cD$ the
eigenspace corresponding to the first non-zero Neumann eigenvalue of $K$.
\begin{enumerate}
\item[(i)] If $K$ is unconditional, then there
exist $i=1,\ldots,n$ and an eigenfunction $0 \not \equiv \vphi \in E_{\lambda_1}$,
such that
$$
 \sigma_i(\vphi) = -\vphi.  $$
\item[(ii)] If $K$ is centrally-symmetric (i.e., $K = -K$), then there
exists an eigenfunction $0 \not \equiv \vphi \in E_{\lambda_1}$,  such that
$$ \vphi(-x) = -\vphi(x) \ \ \ \ \ \text{for} \ x \in K. $$
\end{enumerate}
\label{cor_uncond}
\end{corollary}

\emph{Proof:} Begin with the proof of (i). We are given
the unconditional convex body $K$.
Since $K$ is unconditional, then $f \in E_{\lambda_1}$ implies $\sigma_i(f) \in E_{\lambda_1}$
for $i=1,\ldots,n$. Begin with any non-zero
eigenfunction $f_0 \in E_{\lambda_1}$, and recursively define
$$ f_i = f_{i-1} + \sigma_i(f_{i-1}). $$
Then $f_0,f_1,\ldots,f_n \in E_{\lambda_1}$. If
there exists $i=1,\ldots,n$ such that
$f_i \equiv 0$ then we are done: Suppose $i$ is the minimal
such index. Then $0 \not \equiv f_{i-1} \in E_{\lambda_1}$ with $\sigma_{i-1}(f_{i-1}) = -f_{i-1}$,
and we found our desired eigenfunction.

\medskip It remains to deal with the case where $\psi = f_n$ is a non-zero
eigenfunction. Note that $ \sigma_i(\psi) = \psi$ and hence
\begin{equation}
\sigma_i(\partial^i \psi) = -\partial^i \psi \label{eq_1006}
\end{equation}
for $i=1,\ldots,n$. Therefore,
\begin{equation}
 \int_K \nabla \psi = 0.
\label{eq_1156}
\end{equation} In the proof of Corollary \ref{cor_1154} (the first part,
which did not use the uniform strict convexity) we observed that
(\ref{eq_1156}) implies that
 $\partial^1 \psi,\ldots,\partial^n \psi \in E_{\lambda_1}$. Since
$\int_K |\nabla \psi|^2 > 0$, there exists $i=1,\ldots,n$ with
$\partial^i \psi \not \equiv 0$. We see from (\ref{eq_1006}) that
$\partial^i \psi \in E_{\lambda_1}$ is the eigenfunction
we are looking for. This completes the proof of the first part of the
lemma.

\medskip The proof of the second part is similar. Begin with any $0 \not \equiv
f \in E_{\lambda_1}$ and set $\psi(x) = f(x) + f(-x)$. If $\psi \equiv 0$, then $f$ is an odd
function and we are done. Otherwise, $\psi$ is an even function,
hence $\int_K \nabla \psi = 0$. As before, this implies that
$\partial^1 \psi,\ldots
\partial^n \psi$ are all odd eigenfunctions corresponding to the same
eigenvalue $\lambda_1$. \hfill $\square$

\medskip Corollary \ref{cor_1154} and Corollary
\ref{cor_uncond} seem very much expected.  Notably,
 Nadirashvili \cite{nad} has proved that
in two dimensions, the multiplicity of the first non-zero Neumann
eigenvalue is at most $2$ for any simply-connected domain. Our
simple proof of Corollary \ref{cor_1154} is not applicable in such
generality.  Corollary \ref{cor_1154} is related to the ``hot spots''
problem, see, e.g., Burdzy \cite{burdzy}, Jerison and Nadirashvili
\cite{JN} and references therein. A proof of Corollary \ref{cor_uncond} for the
two-dimensional case -- under much more general assumptions than
convexity -- can be found in \cite[Theorem 4.3]{BB}. However, the
proofs of the two-dimensional results mentioned do not seem to
admit easy generalization to higher dimensions. As observed
by Payne and Weinberger \cite{pw2}, Corollary \ref{cor_uncond}
leads to the following comparison principle:

\begin{corollary} Let $K \subset \RR^n$ be an unconditional convex body with
a $C^{\infty}$-smooth boundary. Assume that $R > 0$ is such that
$$ K \subseteq [-R,R]^n = \{ (x_1,\ldots,x_n) \in \RR^n \, ; \, |x_i| \leq R \ \ \text{for} \ \ i=1,\ldots,n \}. $$
Denote by $\lambda_1 > 0$ the first non-zero Neumann eigenvalue of $K$. Then,
$$ \lambda_1 \geq \frac{\pi^2}{R^2}. $$
Equality holds when $K = [-R, R]^n$, an $n$-dimensional cube. \label{cor_327}
\end{corollary}

\emph{Proof:} A well-known, elementary calculation shows that for any $0 < r \leq R$
and a smooth odd function $\psi:[-r,r] \rightarrow \RR$,
\begin{equation}
\frac{\pi^2}{R^2} \int_{-r}^r \psi^2(x) dx \leq
\frac{\pi^2}{r^2} \int_{-r}^r \psi^2(x) dx \leq \int_{-r}^r \left( \frac{d \psi}{dx} \right)^2 dx.
\label{eq_300}
\end{equation}
According to Corollary \ref{cor_uncond}(i), there exists
an index $1 \leq i \leq n$ and a non-zero eigenfunction $\vphi$ corresponding
to $\lambda_1$ such that $\sigma_{i}(\vphi) = -\vphi$. By Fubini's theorem and (\ref{eq_300}),
$$ \frac{\pi^2}{R^2} \int_K \vphi^2 \leq \int_K |\partial^i \vphi|^2 \leq
\int_K |\nabla \vphi|^2 = \lambda_1 \int_K \vphi^2, $$
hence $\lambda_1 \geq \pi^2 / R^2$. \hfill $\square$

\medskip {\emph Remarks.} \begin{enumerate} \item
Corollary \ref{cor_327} shows that the cube satisfies
a certain domain monotonicity principle for the Neumann Laplacian,
at least in the category of unconditional, convex bodies. The Euclidean ball,
for instance,
does not satisfy a corresponding principle. \item
Suppose $K \subset \RR^n$ is an unconditional convex body. Assume
that $K$ is isotropically normalized, i.e., the random vector $X$
which is distributed uniformly in $K$ is isotropically normalized. Corollary \ref{cor_327}
implies the probably non-optimal bound
\begin{equation}
 \lambda_1(K) \geq c / \log^2 (n+1), \label{eq_341}
 \end{equation}
where $\lambda_1(K) > 0$ is the first non-zero Neumann eigenvalue of $K$,
and $c > 0$ is a universal constant. To establish (\ref{eq_341}), consider 
$$ K^{\prime} = K \cap [-R, R]^n,  \ \ \ \ \ \ \ \text{for} \ R = 50 \log (n+1). $$
Use Corollary \ref{cor_327} to deduce the bound $\lambda_1(K^{\prime}) > c / \log^2 (n+1)$. The body $K^{\prime}$ is a good approximation to the body $K$: It is easily proven that 
$$ Vol(K^{\prime}) \geq \left(1 - \frac{1}{n} \right) Vol_n(K). $$
We may thus apply E. Milman's result \cite[Theorem 1.7]{em},
which builds upon the Sternberg-Zumbrun concavity principle \cite{SZ}, to conclude that
$\lambda_1(K) \geq c \lambda_1(K^{\prime})$
and the bound (\ref{eq_341}) follows. 
See \cite{KLS} for a conjectural better bound, without the logarithmic factor.
\end{enumerate}
\section{Unconditional convex bodies}
\label{sec_uncond}

We begin this section with a corollary to the theorems
of Section \ref{sec_laplace}  and Section \ref{sec_transport}.

\begin{corollary} Let $K \subset \RR^n$ be an
unconditional convex body.
\begin{enumerate}
\item[(i)]
Let $\Psi: K \rightarrow \RR$ be an unconditional, continuous function.
 Then,
$$
Var_K(\Psi) \leq \sum_{i=1}^n
\int_K \left( \Psi (x)- \Psi(\cB_i^+(x)) \right)^2 dx.  $$
\item[(ii)] In particular, suppose $f_1,\ldots,f_n: \RR \rightarrow \RR$ are even, continuous functions.
Denote $\Psi(x_1,\ldots,x_n) = \sum_{i=1}^n f_i(x_i)$. Then,
$$
Var_K(\Psi) \leq  \sum_{i=1}^n
\int_K \sup_{s,t \in J_i(x)} \left( f_i(s) - f_i(t) \right)^2 dx,  $$
where $J_i(x) = [q_i^-(\pi_i(x)), q_i^+(\pi_i(x))] \subset \RR$.
That is, $J_i(x)$ is a symmetric interval about the origin
with the same length as $[\cB_i^-(x), \cB_i^+(x)]$.  
\end{enumerate}
\label{cor_203}
\end{corollary}

\emph{Proof:}  Begin with (i). By approximation, we may assume that $K$ has
a $C^{\infty}$-smooth boundary, and that $\Psi$ is a $C^{\infty}(K)$-smooth
function.
Lemma \ref{main_laplace}
states that
$$  Var_K(\Psi) \leq \sum_{i=1}^n \| \partial^i \Psi \|_{H^{-1}(K)}^2. $$
Fix $i=1,\ldots,n$.
We may apply Theorem \ref{thm_258} for $h = \partial^i \Psi$ since $\int_K \partial^i \Psi = 0$,
as implied by the symmetries of $\Psi$. We may apply
 Lemma \ref{lem_1117},
since clearly $\Psi \left( \cB_i^{+}(x) \right) =
\Psi \left( \cB_i^{-}(x) \right)$ for any $x \in K$.
 Theorem \ref{thm_258} and Lemma \ref{lem_1117} entail the inequality
$$   \| \partial^i \Psi \|_{H^{-1}(K)}^2 \leq \int_K \left( \Psi(x) - \Psi(\cB_i^+(x)) \right)^2 dx. $$
This proves (i). To deduce (ii), denote $\Psi_i(x_1,\ldots,x_n) = f_i(x_i)$. Observe that
$\Psi(x) = \sum_{i=1}^n \Psi_i(x)$ is unconditional and that
for any $x \in K, i=1,\ldots,n$,
$$ \left| \Psi (x)- \Psi(\cB_i^+(x)) \right| = \left| \Psi_i(x) - \Psi_i(\cB_i^+(x)) \right| \leq 
\sup_{s,t \in J_i(x)} \left| f_i(s) - f_i(t) \right|. $$
Thus (ii) follows from (i).
\hfill $\square$

\medskip
We will use the following simple identities:
\begin{equation}
 \int_{-r}^r \left( a |t|^p - a r^p \right)^2 dt
= \frac{2 p^2}{p+1} \int_{-r}^r ( a |t|^p )^2 dt,
\label{eq_217}
\end{equation}
\begin{equation}
 \int_{-r}^r \left( 2  a r^p \right)^2 dt
= 8 a^2 r^{2 p + 1} = 4 (2 p + 1) \int_{-r}^r ( a |t|^p )^2 dt,
\label{eq_333}
\end{equation}
valid for all $a, p, r \geq 0$.

\begin{lemma} Let $X = (X_1,\ldots,X_n)$ be a random vector in $\RR^n$,
that is distributed according to an
unconditional, log-concave density.
Let $p_1,\ldots,p_n > 0$ and let $a_1,\ldots,a_n \geq 0$. Then,
\begin{enumerate}
 \item[(i)] $\displaystyle Var \left( \sum_{i=1}^n a_i |X_i|^{p_i} \right) \leq \sum_{i=1}^n  \frac{2 p_i^2}{p_i + 1}
a_i^2 \, \EE |X_i|^{2p_i}$.
\item[(ii)]
Furthermore, suppose $f_1,\ldots,f_n : \RR \rightarrow \RR$ are even, measurable functions
with $|f_i(t)| \leq a_i |t|^{p_i}$ for all $t \in \RR, \, i=1,\ldots,n$. Then,
$$
 Var \left( \sum_{i=1}^n f_i(X_i) \right) \leq 4 \sum_{i=1}^n  (2 p_i + 1) a_i^2 \, \EE |X_i|^{2p_i}.
 $$
\end{enumerate}
\label{lem_221}
\end{lemma}

\emph{Proof:} Suppose first that $X$ is distributed uniformly in
an unconditional convex body $K \subset \RR^n$.
For $x = (x_1,\ldots,x_n) \in \RR^n$, denote
$$ \Psi(x_1,\ldots,x_n) = \sum_{i=1}^n a_i |x_i|^{p_i}. $$
 The desired bound (i) is equivalent to
$$ Var_K(\Psi) \leq \sum_{i=1}^n
 \frac{2 p_i^2}{p_i + 1}
\int_K a_i^2
 |x_i|^{2p_i} dx_1\ldots dx_n. $$
According to Corollary \ref{cor_203}(i), it suffices to prove that for any $i=1,\ldots,n$,
\begin{equation}
\int_K \left( \Psi(x) - \Psi(\cB_i^+(x)) \right)^2 dx =
 \frac{2 p_i^2}{p_i +1 }
\int_K
a_i^2 |x_i|^{2p_i} dx_1\ldots dx_n.
\label{eq_1147}
\end{equation}
Fix $i=1,\ldots,n$. We will prove (\ref{eq_1147}) by Fubini's
theorem. Fix a point $$ x^{\prime} = (x_1,\ldots,x_{i+1},
x_{i+1},\ldots,x_n) \in \pi_i(K) $$  and denote $r =
q_i^+(x^{\prime}) \geq 0$.
In order to prove (\ref{eq_1147}), it is
enough to show that
$$
\int_{-r}^r \left[
 \sum_{j=1}^n a_j |x_j|^{p_j} - \left( a_i r^{p_i} + \sum_{j\neq i} a_j |x_j|^{p_j}
\right)  \right]^2 dx_i =
 \frac{2 p_i^2}{p_i + 1}
\int_{-r}^r a_i^2 |x_i|^{2p_i} dx_i. $$
The equality we need is exactly the content of (\ref{eq_217}). The proof
of (i) is thus complete, in the case where $X$ is distributed uniformly
in a convex body. The proof of (ii) is almost entirely
identical. By approximation, we may assume that $f_1,\ldots,f_n$ are continuous.
According to Corollary \ref{cor_203}(ii), it is sufficient to prove that
$$
\int_K  \sup_{t, s \in J_i(x)} \left( f_i(s) - f_i(t) \right)^2 dx \leq 4 (2 p_i + 1) \int_K
a_i^2 |x_i|^{2p_i} dx_1\ldots dx_n. $$
This follows by Fubini's theorem and (\ref{eq_333}). The lemma is thus proven,
in the case where $X$ is distributed uniformly in
an unconditional convex body.

\medskip The general case follows via a standard argument.
Let $f: \RR^n \rightarrow [0, \infty)$ stand for the
unconditional, log-concave density
of $X$. Next, we suppose that $f$ is $s$-concave for some integer $s \geq 1$. That is, assume
that
$$ f^{1/s}(\lambda x + (1 - \lambda) y) \geq \lambda f^{1/s}(x) + (1 - \lambda) f^{1/s}(y) $$
for all $0 < \lambda < 1$ and $x,y \in \RR^n$ for which $f(x), f(y) > 0$.
Denote $N = n + s$. For $z \in \RR^N$ we use the coordinates
 $z = (x,y) \in \RR^n \times \RR^s$.
Let $K \subset \RR^N = \RR^n \times \RR^s$ be the unconditional convex body defined by
$$ K = \left \{ (x,y) \ ; \ x \in \RR^n, \, y \in \RR^s, \, |y| \leq \kappa_s ^{-1/s} f^{1/s}(x)  \right \}, $$
where $\kappa_s = \pi^{s/2} / \Gamma(s/2 + 1)$ is the volume of the $s$-dimensional
Euclidean unit ball. Suppose that $Z = (Z_1,\ldots,Z_N)$ is a random vector that is distributed
uniformly in $K$. According to the case already considered,
 conclusions (i) and (ii) hold when the $X_1,\ldots,X_n$
are replaced by $Z_1,\ldots,Z_n$. However, the random vector $(Z_1,\ldots,Z_n)$
has the same distribution as $X = (X_1,\ldots,X_n)$. Thus
(i) and (ii) hold also in the case where the density $f$ is $s$-concave.

\medskip Finally,  an approximation argument eliminates the requirement
that the density of $f$
be $s$-concave: Write $f = e^{-\psi}$ for the unconditional, log-concave density of
$X$. Then, for any $s
> 0$, the function
$$ x \mapsto \left(1 - \frac{\psi(x)}{s} \right)_+^s $$
is unconditional and $s$-concave, where $x_+ = \max \{ x, 0 \}$.
This density clearly tends to $e^{-\psi}$ weakly (and also
uniformly in $\RR^n$) when $s \rightarrow
\infty$. We thus deduce the general case as a limit of the $s$-concave case.\hfill $\square$

\medskip Lemma \ref{lem_221}
may be viewed as a substitute for the sub-independent coordinates
idea of Anttila, Ball and Perissinaki \cite{abp}: Note the absence
of cross terms from the right-hand side of Lemma \ref{lem_221}(i).
Suppose
$X$ is a real-valued random variable  with an even, log-concave
density. A classical inequality (see, e.g., \cite{mp}, or
\cite[Theorem 12]{BaK} and references therein) states that for any
$p \geq 2$,
\begin{equation}
 \left( \frac{\EE |X|^p}{\Gamma(p+1)} \right)^{1/p} \leq
\sqrt{\frac{\EE |X|^2}{2}} \leq  \EE |X|, \label{eq_151}
\end{equation}
where $\Gamma(p+1) = \int_0^{\infty} t^p e^{-t} dt$. For a vector $x = (x_1,\ldots,x_n)
\in \RR^n$ and for $p \geq 1$ we write
$$ \| x \|_p = \left( \sum_{i=1}^n |x_i|^p \right)^{1/p}. $$
The following corollary contains a few obvious consequences of
Lemma \ref{lem_221}.

\begin{corollary} Let $X = (X_1,\ldots,X_n)$ be a random vector in $\RR^n$,
with $\EE X_i^2 = 1$ for $i=1,\ldots,n$,
that is distributed according to an
unconditional, log-concave density. Let $a_1,\ldots,a_n \geq 0$. Then,
\renewcommand{\theequation}{i}
\begin{equation}
 Var \left( \sum_{i=1}^n a_i X_i^2 \right) \leq C^{\prime} \sum_{i=1}^n a_i^2,
\end{equation}
where $C^{\prime} \leq 16$ is a universal constant. Consequently,
\renewcommand{\theequation}{ii}
\begin{equation}
Var( |X|^2 ) \leq C^2 n \ \ \  \ \ \ \ \text{and} \ \ \ \  \ \ \
\EE \left( |X| - \sqrt{n} \right)^2 \leq C^2,
\end{equation}
with $C \leq 4$, a positive universal constant. Moreover, for any
$p \geq 1$,
\renewcommand{\theequation}{iii}
\begin{equation}
 \sqrt{Var \left( \| X \|_p \right)} \leq C_p n^{\frac{1}{p} - \frac{1}{2}}
\end{equation}
where $C_p > 0$ is a constant depending only on $p$.
\label{cor_204}
\end{corollary}
\renewcommand{\theequation}{\arabic{equation}}

\emph{Proof:}
According to the Pr\'ekopa-Leindler inequality (see, e.g., the
first pages of \cite{pisier_book}), the random variable $X_i$ has
an even, log-concave density for all $i$. From Lemma \ref{lem_221}(i) and
(\ref{eq_151}) we see that
\begin{equation} \label{eq_441}
  Var \left( \sum_{i=1}^n a_i X_i^2
 \right)  \leq  \frac{8}{3} \sum_{i=1}^n a_i^2 \EE |X_i|^4 \leq 16  \sum_{i=1}^n a_i^2 \left( \EE |X_i|^2 \right)^2 = 16 \sum_{i=1}^n a_i^2. \nonumber
\end{equation}
This proves (i). By setting $a_i = 1 \ (i=1,\ldots,n)$
in (\ref{eq_441}), we deduce that
$$ \EE \left( |X| - \sqrt{n} \right)^2 \leq \frac{1}{n}
\EE \left( |X| - \sqrt{n} \right)^2 \cdot \left( |X| + \sqrt{n}
\right)^2 = \frac{1}{n} \EE \left( |X|^2 - n \right)^2 \leq 16, $$
and (ii) is proven. Denote $E =  \EE \| X \|_p^p$.
From Lemma \ref{lem_221}(i) and (\ref{eq_151}) we conclude that
$$ \EE \left( \| X \|_p^p - E \right)^2
= Var \left( \sum_{i=1}^n |X_i|^p \right) \leq 2^{1 - p} p
\Gamma(2 p + 1) n. $$
For any $p \geq 2$, we have $\EE |X_i|^p \geq (\EE X_i^2)^{p/2} = 1$. For $1 \leq p \leq 2$,
$$ \EE |X_i|^p \geq \left( \EE |X_i| \right)^p \geq
2^{-p/2} \left( \EE X_i^2 \right)^{p/2} = 2^{-p/2} \geq 2^{-1/2}, $$ according
to (\ref{eq_151}). Hence, $E = \sum_i \EE |X_i|^p \geq n / \sqrt{2}$ and
$$ Var \left( \| X \|_p \right) \leq \EE \left( \| X \|_p - E^{1/p} \right)^2
\leq E^{-2\frac{p-1}{p}} \EE \left( \| X \|_p^p - E \right)^2 \leq
C_p n^{2/p-1}, $$ where $C_p$ is a constant depending solely on $p
\geq 1$.  This completes the proof. \hfill $\square$

\medskip Schechtman and Zinn \cite{SZ1, SZ2} provided estimates related to
 Corollary \ref{cor_204} for the case where $X$ is distributed uniformly
in the unit ball $\{ x \in \RR^n ; \| x \|_q \leq 1 \}$, for $q
\geq 1$. More information regarding unconditional, log-concave
densities in high dimension, especially in the large deviations
scale, is available from Bobkov and Nazarov \cite{BN1, BN2}. Under
the assumptions of Corollary \ref{cor_204}, they showed, for
instance, that
$$ \PP \left( \frac{1}{\sqrt{n}} \sum_{i=1}^n X_i \geq t \right) \leq C \exp \left(-c t^2 \right)
 \ \ \ \ \ \ \forall t \geq 0, $$
where $c, C > 0$ are universal constants. Another large-deviations
estimate that was proved by Bobkov and Nazarov \cite{BN1, BN2} is
that
\begin{equation}
 \PP \left( |X| \geq t \right) \leq C \exp \left(-c t \right)
\ \ \ \ \ \ \text{for} \ t \geq C \sqrt{n}.
\label{eq_448}
\end{equation}
 Paouris \cite{pa_cr, pa_gafa} was remarkably able to
generalize  inequality (\ref{eq_448}) to the class of all
isotropically-normalized random vectors with a log-concave density in $\RR^n$. Regarding smaller
values of $t$ in (\ref{eq_448}), the currently known bounds, which are valid for all isotropically-normalized, log-concave random vectors, are of
the form
\begin{equation} \PP \left( \left| \frac{|X|}{\sqrt{n}} - 1 \right| \geq t \right) \leq C \exp
 \left(-c n^{\alpha} t^{\beta} \right)
\ \ \ \ \ \ \text{for} \ 0 < t < 1, \label{eq_454} \end{equation}
with, say, $\alpha = 0.33$ and $\beta = 3.33$ (see \cite{power_law}).

\medskip Cordero-Erausquin, Fradelizi and Maurey \cite{cfm} have recently
proved the so-called (B)-conjecture in the unconditional
case. This entails the following improvement over the
Brunn-Minkowski theory:
\begin{itemize}
\item The function $ t \mapsto \PP \left( |X| \leq e^t \right)$ is log-concave in $t \in \RR$.
\end{itemize}
(The Pr\'ekopa-Leindler inequality leads to the weaker statement
in which the $e^t$ is replaced by $t$).
Corollary \ref{cor_204}(ii) and Markov-Chebychev's inequality yield
$$
 \PP \left(  |X|  \leq \sqrt{n} - 8 \right) \leq
\frac{1}{4}, \ \ \ \ \PP \left(  |X| \leq \sqrt{n}  + 8
 \right) \geq \frac{3}{4}. $$ 
The log-concavity of the map
$ s \mapsto \PP \left(  |X| \leq e^s \right)$
thus implies that for any $t \geq 0$,
$$ \PP \left(  |X| \leq  (\sqrt{n} - 8 )
\cdot \left( \frac{\sqrt{n} - 8 }{\sqrt{n} + 8 } \right)^t \right)
\leq \frac{1}{4 \cdot 3^t}. $$ After some simple manipulations, we deduce the inequality \begin{equation}
\PP \left(  |X| \leq \sqrt{n} - t \right) \leq C \left( 1
 - \frac{t}{\sqrt{n}} \right)^{c \sqrt{n}} \leq C \exp(-c t),
\label{eq_731}
\end{equation}
valid for all $0 \leq t \leq \sqrt{n}$, for some universal constants $c, C > 0$. We currently do not know how to prove a bound
as in (\ref{eq_731}) for the probability $\PP \left( |X| \geq
\sqrt{n} + t \right)$. The weaker estimate 
$$ \PP \left( |X| \geq \sqrt{n} + t \right) \leq C \exp \left( -c \sqrt{t} \right) $$
follows by combining Corollary \ref{cor_204}(ii) with the distribution inequalities of Nazarov, Sodin and Volberg \cite{nsv}. We omit the details.

\section{Berry-Esseen type bounds}
\label{sec_berry}

In previous sections we established sharp {\it thin shell estimates}
for unconditional, log-concave densities. In the present section
we complete the proof of Theorem \ref{thm_511}. 
The argument we present is quite technical
and is very much related to classical
treatments of the central limit theorem for independent random variables.
The reader may refer to, e.g., \cite[Vol. II, Chapter XVI]{feller} for background on the
rate of convergence in the classical central limit
theorem.
We are indebted to Sasha Sodin for many discussions, suggestions
and simplifications that have lead to the proofs we present below.

\medskip Before proceeding to the actual proof, let us describe the general idea. 
Introduce independent, symmetric Bernoulli variables $\Delta_1,\ldots,
\Delta_n$. That is,
$$ \PP (\Delta_i = 1) = \PP (\Delta_i = -1) = 1/2 \ \ \ \ \ \ \ (i=1,\ldots,n). $$
These Bernoulli variables are also assumed to be independent of $X$. 
Write
$$ \vphi(t) = \frac{1}{\sqrt{2 \pi}} e^{-t^2/2}
\ \ \  \ \ \ \ \ \text{and} \ \ \ \ \ \ \ \ \ \ \ \Phi(t) = \int_t^{\infty} \vphi(s) ds $$
for all $t \in \RR$. We condition on $X$, and apply the classical Berry-Esseen bound to obtain 
\begin{equation}
 \left| \PP  \left( \frac{\sum_i \Delta_i X_i}{\sqrt{n}} \geq t \right)
- \Phi \left(t \sqrt{n} / |X| \right) \right| \leq C \frac{\sum_i |X_i|^3}{\left( \sum_i |X_i|^2 \right)^{3/2}} \leq \frac{C^{\prime}}{\sqrt{n}}  \label{eq_1111_}
\end{equation}
where the last inequality holds only for ``typical'' values of $X$. Since $|X| / \sqrt{n}$ is strongly concentrated around $1$, as we learn from (\ref{eq_429_}), we may substitute the 
$\Phi \left(t \sqrt{n} / |X| \right)$ term in (\ref{eq_1111_}) by $\Phi(t)$.
Observe that since $X$ is unconditional, the random variables
$$ \sum_i X_i \ \ \ \ \ \ \ \text{and} \ \ \ \ \ \ \ \sum_i X_i \Delta_i $$
have exactly the same distribution. Hence, by considering the expectation over $X$ in (\ref{eq_1111_}), we deduce a weaker 
version of (\ref{eq_1124}) where the $C / n$ is replaced with $C / \sqrt{n}$. 
In order to arrive at the optimal bound, we need to apply a smoothing technique: The estimate (\ref{eq_1111_}) will be replaced with a much better Berry-Esseen inequality which is available for the random variable
$\Gamma + \left( \sum_i \Delta_i X_i \right) \left / \sqrt{n} \right.$, for an appropriate ``small'' random variable $\Gamma$. The details will be described next.

\medskip Throughout this section, we fix a symmetric random variable $\Gamma$ with $\EE \Gamma^6 < \infty$,
independent of everything else, such that the even function $\gamma(\xi) =
\EE \exp(-i \xi \Gamma)$ satisfies
\begin{equation}
\gamma(\xi) = 0 \ \ \ \ \ \text{for} \ \ \ \ \ |\xi| \geq 1 \label{eq_1032}
\end{equation}
and
\begin{equation}
1 - 1000 \xi^2 \leq \gamma(\xi) \leq 1 \ \ \ \ \ \text{for} \ \ \ \ \ \xi \in \RR. \label{eq_313}
\end{equation}
For instance, $\Gamma$ may be the random variable whose density is $$ x \mapsto \kappa_1 \sin^{8}(\kappa_2 x) / x^{8}, $$ for appropriate universal constants $\kappa_1, \kappa_2$. (For this specific choice, $\gamma$ is the $8$-fold
convolution of the characteristic function of an interval.)
We shall use the standard $O$-notation in this section. The notation $O(x)$, for some expression $x$,
is an abbreviation for some complicated quantity $y$ with the property that
$$ |y| \leq C x $$
for some universal constant $C > 0$. All constants hidden in the $O$-notation in our
proof are in principle explicit. 
The following lemma seems rather standard (see \cite[Vol. II, Chapter XVI]{feller} for similar statements).
For lack of a precise reference, we provide its proof.

\begin{lemma} Suppose $\Delta_1,\ldots,\Delta_n$ are independent, symmetric Bernoulli random variables. 
Let $0 \neq \theta = (\theta_1,\ldots,\theta_n) \in \RR^n$ and $\sigma > 0$. Assume that
\begin{equation}
 \sum_{i; |\theta_i| \geq \sigma} \theta_i^2 \leq \frac{1}{2} |\theta|^2. \label{eq_955} \end{equation}
Then, for any $t \in \RR$,
\begin{equation}
 \left|   \PP \left( \sigma \Gamma + \sum_{i=1}^n \theta_i \Delta_i \geq t \right) \,  - \, \Phi \left( \frac{t}{|\theta|} \right) \right| \leq C \left(   \frac{\sigma^2}{|\theta|^2}  + \sum_{i=1}^n \frac{\theta_i^4}{|\theta|^4} \right), \label{eq_1022}
\end{equation}
where $C > 0$ is a universal constant. \label{lem_700}
\end{lemma}

\emph{Remark.} Note that when $\theta_i = 1/\sqrt{n} = \sigma$ for all $i$, the error term in Lemma \ref{lem_700} is $O(1/n)$. The addition of $\Gamma / \sqrt{n}$ allows us to deduce a better bound than the $O(1/\sqrt{n})$ guaranteed by the Berry-Esseen inequality.

\medskip
\emph{Proof of Lemma \ref{lem_700}:} The validity of both the assumptions and the conclusions
of the lemma is not altered if we replace $\theta, \sigma$ with $r \theta, r \sigma$ for any $r > 0$.
Normalizing, we may assume that $|\theta| = 1$. By symmetry, it is enough to prove (\ref{eq_1022}) for non-negative $t$.
 Fix $t \geq 0$.
Observe that for any $\xi \in \RR$, $$ \EE \exp \left(-i \xi \left[\sigma \Gamma + \sum_{i=1}^n \theta_i \Delta_i \right]\right) = \gamma \left(  \sigma \xi \right) \prod_{i=1}^n \cos(\theta_i \xi). $$
Thus, from the Fourier inversion formula
(see, e.g., \cite[Vol. II, Chapter XVI]{feller}),
\begin{eqnarray} \nonumber
\lefteqn{ \PP \left(  \sigma \Gamma + \sum_{i=1}^n \theta_i \Delta_i \leq t \right) \,  - \, \frac{1}{\sqrt{2 \pi}} \int_{-\infty}^t \exp(-s^2 / 2) ds } \\ & = &
\frac{1}{2 \pi} \int_{-\infty}^{\infty} \left[  \gamma \left(  \sigma \xi  \right) \prod_{i=1}^n \cos(\theta_i \xi) - e^{-\xi^2 / 2} \right] \frac{e^{i t \xi} - 1}{i \xi}  d \xi. \label{eq_1028}
\end{eqnarray}
Denote $\eps = \sqrt{\sum_i \theta_i^4}$.
To prove the lemma, it suffices to bound the absolute value of the integral in (\ref{eq_1028}) by $C^{\prime} (\eps^2 + \sigma^2)$. We express the integral in (\ref{eq_1028}) as
$I_1 + I_2 + I_3$ where $I_1$ is the integral over $\xi \in [-\eps^{-1/2},\eps^{-1/2}]$,
 $I_2$ is the integral over $\eps^{-1/2} \leq |\xi| \leq \sigma^{-1}$
(when $\eps^{-1/2} > \sigma^{-1}$, we set $I_2 = 0$)
and $I_3$ is the integral over  $|\xi| \geq \max \{\sigma^{-1}, \eps^{-1/2}\}$.

\medskip Begin with estimating $I_1$. We use the elementary inequality
$$ e^{s^2/2} \cos s = e^{O( s^4)} \ \ \ \ \ \text{for} \ \ \  |s| \leq 1. $$
Since $|\theta_i| \leq \eps^{1/2}$ for all $i$, then for $|\xi| \leq \eps^{-1/2}$,
\begin{equation}
 \left| \prod_{i=1}^n e^{\xi^2 \theta_i^2 / 2} \cos(\theta_i \xi) - 1 \right| =
 \left| e^{ O \left(\xi^4 \sum_{i=1}^n  \theta_i^4 \right) } - 1 \right| \leq C^{\prime} \xi^4 \eps^2. \label{eq_311}
\end{equation}
 Combine (\ref{eq_311}) with (\ref{eq_313}) to deduce that for $|\xi| \leq \eps^{-1/2}$,
$$  \gamma \left(  \sigma \xi \right) \prod_{i=1}^n e^{\xi^2 \theta_i^2 / 2} \cos(\theta_i \xi) = \left(1 + O (  \sigma^2 \xi^2 ) \right) \left(1 + O(  \xi^4 \eps^2) \right) = 1 +  O(\sigma^2 \xi^2 + \xi^4 \eps^2). $$
The latter estimate yields
\begin{eqnarray} \nonumber
|I_1| & = & \left|  \int_{-\eps^{-1/2}}^{\eps^{-1/2}} e^{-\xi^2 / 2} \left[ \gamma \left(  \sigma \xi \right) \prod_{i=1}^n e^{\xi^2 \theta_i^2 / 2} \cos(\theta_i \xi) - 1 \right] \frac{e^{i t \xi} - 1}{i \xi} d \xi \right|  \\ \nonumber  & \leq & C^{\prime}
\int_{-\infty}^{\infty} e^{-\xi^2 / 2} \left( \sigma^2 \xi^2 + \xi^4 \eps^2 \right) \frac{2}{|\xi|} d \xi \, \, \leq  \, \, \tilde{C}
\left( \sigma^2 + \eps^2 \right), \label{eq_343} \end{eqnarray}
since $0 < \eps \leq 1$.

\medskip Next we estimate $I_2$, in the case where $\eps^{-1/2} \leq \sigma^{-1}$ (in the complementary
case, $I_2 = 0$). Denote $ \cI = \left \{ 1 \leq i \leq n \, ; \, |\theta_i| \leq \sigma \right \}$. Then, by (\ref{eq_955}),
\begin{equation} \sum_{i \in \cI} \theta_i^2 \geq 1/2.
\label{eq_401}
\end{equation}
We will use the elementary inequality $|\cos s| \leq e^{-c s^2}$ for $|s| \leq 1$.
According to (\ref{eq_401}), whenever $|\xi| \leq \sigma^{-1}$,
$$ \left| \prod_{i=1}^n \cos(\theta_i \xi) \right| \leq  \prod_{i \in \cI} |\cos(\theta_i \xi)|  \leq  e^{-c \xi^2 \sum_{i \in \cI} \theta_i^2} \leq e^{-c \xi^2  /2}. $$
 Apply the well-known bound $\int_s^{\infty} e^{-u^2 / 2} \leq C e^{-c s^2}$ for $s \geq 0$, to deduce
\begin{eqnarray}
\nonumber
 |I_2| & \leq  & 2 \int_{\eps^{-1/2}}^{\sigma^{-1}} \left[ \left| \prod_{i=1}^n \cos(\theta_i \xi) \right| + e^{-\xi^2/2} \right] \frac{2}{|\xi|} d \xi  \\ & \leq & 4 \int_{\eps^{-1/2}}^{\sigma^{-1}} \left[ e^{-c \xi^2  /2 } + e^{-\xi^2/2} \right] d\xi \leq \bar{C} e^{-\tilde{c}/\eps} \leq \tilde{C} \eps^2.\label{eq_540_}
\end{eqnarray}
The bound for $I_3$ is easy. From (\ref{eq_1032}) we have $\gamma( \sigma \xi) = 0$ for $|\xi| \geq \sigma^{-1}$. Hence,
$$ |I_3| \leq 2 \left| \int_{\max \{\sigma^{-1}, \eps^{-1/2}\}}^{\infty} e^{-\xi^2/2} \frac{2}{|\xi|} d\xi \right| \leq C e^{-c/\sigma^2} \leq \bar{C} \sigma^2. $$
The lemma follows by combining the above bound for $|I_3|$ with the bound (\ref{eq_540_}) for $|I_2|$ and the bound (\ref{eq_343}) for $|I_1|$. \hfill $\square$

\begin{lemma} Let $X = (X_1,\ldots,X_n)$ be a random vector in $\RR^n$,
with $\EE X_i^2 = 1$ for $i=1,\ldots,n$,
that is distributed according to an
unconditional, log-concave density. Let $(\theta_1,\ldots, \theta_n) \in S^{n-1}$
and denote $\eps = 10 \sqrt{\sum_i \theta_i^4}$. Then,
$$ \displaystyle \PP \left(
 \frac{1}{2} \leq \sum_{i=1}^n \theta_i^2 X_i^2 \leq \frac{3}{2} \ \ \ \ \ \ \ \text{and}
 \ \ \ \sum_{i ; |\theta_i X_i| \geq \eps} \theta_i^2 X_i^2 \leq \frac{1}{4}  \right) \geq 1 - C \eps^2, $$
where $C > 0$ is a universal constant.
\label{lem_426}
 \end{lemma}

\emph{Proof:} Note that $\EE \sum_{i=1}^n \theta_i^2 X_i^2 = 1$. According to the Chebyshev's inequality and  Corollary \ref{cor_204},
\begin{equation}
 \PP \left( \left| \sum_{i=1}^n \theta_i^2 X_i^2 - 1 \right| \geq 1/2 \right)
\leq 4 Var \left( \sum_{i=1}^n \theta_i^2 X_i^2 \right) \leq 64 \sum_{i=1}^n \theta_i^4  \leq \eps^2. \label{eq_340}
\end{equation} Denote $Y = \sum_{i ; |\theta_i X_i| \geq \eps} \theta_i^2 X_i^2$. Clearly,
$$ \eps^2 Y  = \eps^2 \sum_{i ; |\theta_i X_i| \geq \eps} \theta_i^2 X_i^2 \leq \sum_{i=1}^n \theta_i^4 X_i^4. $$
Therefore
$$ \EE Y \leq \eps^{-2} \sum_{i=1}^n \theta_i^4 \EE X_i^4 \leq 6 \eps^{-2} \sum_{i=1}^n \theta_i^4 \leq \frac{1}{10}, $$
where we used the inequality
$\EE X_i^4 \leq 6 (\EE X_i^2)^2 = 6$, quoted above as (\ref{eq_151}).
Next, apply  Lemma \ref{lem_221}(ii) with $f_i(t) = \theta_i^2 t^2$
for $|t| \geq \eps / \theta_i$ and $f_i(t) = 0$ otherwise. According 
to the conclusion of that lemma,$$
Var(Y) = Var \left(  \sum_{i ; |\theta_i X_i| \geq \eps} \theta_i^2 X_i^2 \right)
\leq 4 \sum_{i=1}^n 5 \theta_i^4 \EE X_i^4 \leq 120 \sum_{i=1}^n \theta_i^4 \leq C \eps^2.
$$
Denote $\mu = \EE Y \leq 1/10$.
Another
application of the Chebyshev inequality yields
\begin{equation}
 \PP \left( Y \geq \frac{1}{4} \right) \leq \PP \left( |Y - \mu| \geq \frac{1}{10}
\right) \leq 100 Var(Y) \leq C \eps^2.
\label{eq_429}
\end{equation}
The lemma follows from (\ref{eq_340}) and (\ref{eq_429}).
\hfill $\square$

\begin{lemma}
Let $X = (X_1,\ldots,X_n)$ be a random vector in $\RR^n$,
with $\EE X_i^2 = 1$ for $i=1,\ldots,n$,
that is distributed according to an
unconditional, log-concave density. Let $(\theta_1,\ldots, \theta_n) \in S^{n-1}$
and denote $\eps = 10 \sqrt{\sum_i \theta_i^4}$. Then, for any $t \in \RR$,
$$ \left|   \PP \left(  \eps \Gamma + \sum_{i=1}^n \theta_i X_i \geq t \right) \,  - \, \Phi(t) \right| \leq C \eps^2, $$
where $C > 0$ is a universal constant.
\label{lem_612}
\end{lemma}

\emph{Proof:} We may assume that $\eps$ is smaller than some given positive universal constant,
as otherwise the conclusion is trivial. Let $\Delta_1,\ldots,\Delta_n$ be  independent, symmetric, Bernoulli random variables, that are independent also of $X$. For $t \in \RR$ and $x = (x_1,\ldots,x_n) \in \RR^n$ define 
$$ P(t,x) = \PP \left(  \eps \Gamma + \sum_{i=1}^n \theta_i x_i \Delta_i \geq t \right). $$
Since the density of $X$ is unconditional, the random variable $\sum_i \theta_i X_i$ has the same distribution as $\sum_i \theta_i X_i \Delta_i$. Fix $t \in \RR$. Then,
\begin{equation}
\PP \left( \eps \Gamma + \sum_{i=1}^n \theta_i X_i \geq t \right)
= \PP \left(  \eps \Gamma + \sum_{i=1}^n \theta_i X_i \Delta_i \geq t \right) =  \EE P(t,X).
\label{eq_1044}
\end{equation}
Write $\cA \subset \RR^n$ for the collection of all $x = (x_1,\ldots,x_n) \in \RR^n$
for which
$$
 \frac{1}{2} \leq \sum_{i=1}^n \theta_i^2 x_i^2 \leq \frac{3}{2} \ \ \ \ \ \ \  \text{and}
 \ \ \ \ \ \ \ \ \sum_{i ; |\theta_i x_i| \geq \eps} \theta_i^2 x_i^2 \leq \frac{1}{4} \leq \frac{1}{2} \sum_{i=1}^n \theta_i^2 x_i^2.
$$
We may  apply Lemma \ref{lem_700} for  $(\theta_1 x_1,\ldots, \theta_n x_n)$ and
for $\sigma = \eps$, and conclude that,
$$
\left| \, P(t,x) \, - \, \Phi \left( \frac{t}{\sqrt{\sum_{i=1}^n \theta_i^2 x_i^2}} \right) \, \right| \leq C \left( \eps^2 + \sum_{i=1}^n \theta_i^4 x_i^4 \right)
 \ \ \ \ \ \ \ \ \text{for all} \ x \in \cA.
$$
From Lemma \ref{lem_426} we have $\PP(X \not \in \cA) \leq C \eps^2$.
Consequently,
\begin{eqnarray}
\label{eq_10441}
 \lefteqn{\left| \EE P(t,X) - \EE \Phi \left( \frac{t}{\sqrt{\sum_{i=1}^n \theta_i^2 X_i^2}} \right)
\right|} \\ & \leq & 2 \PP(X \not \in \cA) + C \EE \left( \eps^2 + \sum_{i=1}^n \theta_i^4 X_i^4 \right)
\leq C^{\prime} \eps^2, \nonumber
\end{eqnarray}
where we used once more the bound $\EE X_i^4 \leq 6 (\EE X_i^2)^2 = 6$.
According to (\ref{eq_1044}) and (\ref{eq_10441}),
in order to prove the lemma, all we need is  to show that
\begin{equation}
 \left| \, \EE \Phi \left( \frac{t}{\sqrt{\sum_{i=1}^n \theta_i^2 X_i^2}} \right) \, - \, \Phi(t) \, \right| \leq C \eps^2. \label{eq_507} \end{equation}
Write $Y = \sum_{i=1}^n \theta_i^2 X_i^2$.
Then $\PP(Y \geq 1/2) \geq 1- C \eps^2$, by Lemma \ref{lem_426}.
Therefore, to prove (\ref{eq_507})
and complete the proof of the lemma, it suffices to show that
\begin{equation}
 \EE \left[ \left. \Phi \left( \frac{t}{\sqrt{Y}} \right)  -  \Phi(t) \right| Y \geq 1/2 \right]
= O (\eps^2). \label{eq_5077} \end{equation}
We may assume that $\eps$ does not exceed a small positive universal constant, hence $\PP(Y \geq 1/2)^{-1} \leq (1 - C \eps^2)^{-1} \leq 1 + C^{\prime} \eps^2$. Therefore,
\begin{equation}
 1 = \EE Y \leq \EE \left(  Y \left| Y \geq \frac{1}{2} \right. \right) \leq \PP(Y \geq 1/2)^{-1} \leq 1 + C^{\prime} \eps^2.
\label{eq_5071}
\end{equation}
Corollary \ref{cor_204}(i) implies that $\EE (Y-1)^2 \leq C \eps^2$. Hence,
\begin{equation} \EE \left( \left. \left( Y  - 1 \right)^2 \right| Y \geq \frac{1}{2} \right) \leq
 \EE \left( Y  - 1 \right)^2 / \PP(Y \geq 1/2) \leq \tilde{C} \eps^2.
\label{eq_508}
\end{equation}
Denote $F(u) = \Phi(t / \sqrt{u})$. Clearly, $\vphi(s) s = O(1)$ and $\vphi^{\prime}(s) s^2 = O(1)$
for any $s \in \RR$. Consequently, for any $u \geq 1/2$,
$$ F^{\prime}(u) = \frac{1}{2 u} \vphi \left( \frac{t}{\sqrt{u}} \right) \frac{t}{\sqrt{u}}  = O(1) $$
and
$$ F^{\prime \prime}(u) = -\frac{3}{4 u^2}
\vphi \left( \frac{t}{\sqrt{u}} \right) \frac{t}{\sqrt{u}} - \frac{1}{4 u^2} \vphi^{\prime} \left( \frac{t}{\sqrt{u}} \right) \frac{t^2}{u} = O(1). $$
By Taylor's theorem,
\begin{eqnarray*}
\lefteqn{ \EE \left[ \left. \Phi \left( t / \sqrt{Y} \right)  -  \Phi(t) \right| Y \geq 1/2 \right]
\, \, = \, \, \, \EE \left[ F(Y) - F(1) \left| Y \geq 1/2 \right. \right]} \\ & = & \EE \left[ \left. F^{\prime}(1) (Y - 1) + O \left( (Y - 1)^2 \right) \right| Y \geq 1/2 \right]
\\ & = & F^{\prime}(1) \left( \EE (Y-1) \left | Y \geq \frac{1}{2} \right. \right)  + O(\eps^2)
= O(\eps^2), \phantom{aaaaaaaaa} \end{eqnarray*}
where we used the estimates for $F^{\prime}, F^{\prime \prime}$
and the bounds (\ref{eq_5071}) and (\ref{eq_508}). This completes
the proof of (\ref{eq_5077}). The lemma is proven. \hfill $\square$

\medskip
Our next goal is to eliminate the ``$\eps \Gamma$'' term from the conclusion
of Lemma \ref{lem_612}. The following short computational lemma serves this purpose. We shall use the standard estimate
\begin{equation}
 c \frac{\vphi(t_0)}{t_0+1}  \leq \Phi(t_0) \leq C \frac{\vphi(t_0)}{t_0+1} \leq \bar{C} \vphi(t_0) \label{eq_640}
 \end{equation}
 for any $t_0 \geq 0$ (see, e.g., \cite[Vol. I, Section VII.1]{feller}).

\begin{lemma} Let $t_0 \geq 0$ and denote $\delta = \Phi(t_0)$. Then,
\begin{enumerate}
\item[(i)] $\displaystyle \Phi \left(t_0 + 2 \delta^{1/4} \right) \geq C_1^{-1} \delta$.
\item[(ii)] $\displaystyle 1 - \Phi \left(t_0 - 2 \delta^{1/4} \right) \geq 1 - \Phi(-2) \geq C_1^{-1} \geq C_1^{-1} \delta$.
 \item[(iii)] Suppose $x > 0$ satisfies $\displaystyle \left|\frac{1}{x} - \frac{1}{\vphi(t_0)} \right| \leq c_2 \delta^{-3/4}$.
 Then $\displaystyle x^2 \leq C_1 \delta$.
\end{enumerate} \label{lem_1034}
Here, $C_1 > 1$ and $0 < c_2 < 1$ are universal constants.
 \end{lemma}

\emph{Proof:} We have $t_0 \delta^{1/4} \leq C t_0 (\vphi(t_0))^{1/4} \leq C^{\prime}$
according to (\ref{eq_640}). Hence,
$$ \frac{\Phi\left(t_0 + 2 \delta^{1/4} \right)}{\Phi(t_0)} \geq c^{\prime} \exp \left[ \frac{t_0^2}{2} -
\frac{\left(t_0 + 2 \delta^{1/4} \right)^2}{2} \right] \geq \hat{c} \exp \left( -2 t_0 \delta^{1/4} \right) \geq c^{\prime}, $$
and (i) is proven. The statement (ii) is self-explanatory. Regarding (iii),
it is readily verified that $\tilde{c} (t_0 + 1)^{3/4} \leq \vphi(t_0)^{-1/4}$ for any $t_0 \geq 0$. Therefore, by (\ref{eq_640}), for a sufficiently small $c_2 > 0$,
$$ \frac{1}{\vphi(t_0)} - \frac{c_2}{\delta^{3/4}} \geq \frac{1}{\vphi(t_0)} - \frac{\tilde{c}(t_0+1)^{3/4}}{2\vphi(t_0)^{3/4}} \geq
\frac{1}{\vphi(t_0)} - \frac{\vphi(t_0)^{-1/4}}{2 \vphi(t_0)^{3/4}} = \frac{1}{2 \vphi(t_0)}.  $$
Note also that $\vphi(t_0) \leq C / (t_0+1)$. Consequently, for any $x > 0$,
$$ \left| \frac{1}{x} - \frac{1}{\vphi(t_0)} \right| \leq \frac{c_2}{\delta^{3/4}}
\ \ \ \ \ \Rightarrow  \ \ \ \ \ x \leq 2 \vphi(t_0) \leq C \sqrt{\frac{\vphi(t_0)}{t_0 + 1}} \leq \tilde{C} \sqrt{\delta}, $$
where we used (\ref{eq_640}) again.
\hfill $\square$

\begin{lemma} Let $X$ be a real-valued random variable with
an even, log-concave density. Let $0 < \eps < 1, A \geq 1$. Suppose that
for any $t \in \RR$,
\begin{equation} \left|   \PP \left( \eps \Gamma + X  \geq t \right)
\, - \, \Phi(t)  \right| \leq A \eps^2. \label{eq_627} \end{equation}
Then, for any $t \in \RR$,
\begin{equation}
 \left|   \PP \left( X  \geq t \right)
\, - \, \Phi(t) \right| \leq C A \eps^2,  \label{eq_520}
\end{equation}
where $C > 0$ is a universal constant.
\label{lem_618}
\end{lemma}

\emph{Proof:} By approximation, we may assume that the density of $X$ is $C^1$-smooth
and everywhere positive (e.g., convolve $X$ with a very small gaussian).
We may also assume that $\eps \leq c$ for a small universal constant $c > 0$. The function
$$E(t) = \left|   \PP \left( X  \geq t \right)
\, - \, \Phi(t) \right| \ \ \ \ \ \ \ \ \ \ \ (t \in \RR) $$
is continuous and vanishes at $\pm \infty$. Consequently, there exists
$t_0 \in \RR$ where $E(t)$ attains its maximum. Since $E$ is an even function,
we may assume that $t_0 \geq 0$.
Write $f: \RR \rightarrow [0, \infty)$ for the density of $X$. As $E^{\prime}(t_0) = 0$,
\begin{equation}
 f(t_0) = \vphi(t_0) = \frac{1}{\sqrt{2 \pi}} e^{-t_0^2 / 2}.
\label{eq_114}
\end{equation}
To prove the lemma, it suffices to show that $\max_t E(t) = E(t_0) \leq C A \eps^2$.

\medskip
{\bf Step 1:} Suppose first that $\Phi(t_0) \leq 2 C_1 A \eps^2$, for $C_1$ being the universal
constant from Lemma \ref{lem_1034}. Then by (\ref{eq_627}),
$$ \PP \left( \eps \Gamma + X \geq t_0 \right) \leq \Phi(t_0) + A \eps^2 \leq (2 C_1 + 1) A \eps^2, $$
hence,
$$ \PP(X \geq t_0) = 2 \PP (X \geq t_0, \Gamma \geq 0) \leq 2 \PP \left( \eps \Gamma + X \geq t_0 \right) \leq (4 C_1 + 2) A \eps^2. $$
Consequently, since $\Phi(t_0) \leq 2 C_1 A \eps^2$,
$$ \max_{t \in \RR} E(t) = E(t_0) = \left| \PP \left(  X \geq t_0 \right) - \Phi(t_0) \right| \leq (6 C_1 + 2) A \eps^2 \leq \bar{C} A \eps^2. $$
The desired estimate (\ref{eq_520}) is therefore proven, in the case where $\Phi( t_0 ) \leq 2 C_1 A \eps^2$.

\medskip
{\bf Step 2:}  It remains to deal with the case where $t_0 \geq 0$ satisfies $\Phi(t_0) > 2 C_1 A \eps^2$. Denote $\delta = \Phi(t_0) \geq 2 C_1 A \eps^2 \geq A \eps^2$. Note that
\begin{equation} \PP \left(
|\eps \Gamma| \geq \delta^{1/4} \right) \, \leq \,  \frac{\eps^6 \EE \Gamma^6}{\left(\delta^{1/4} \right)^6} \leq  C \frac{\eps^{3}}{A^{3/2}} \leq C \eps \delta \leq \frac{\delta}{4 C_1} \label{eq_1116}
\end{equation}
under the legitimate assumption that $\eps$ is smaller than a given universal constant.
From Lemma \ref{lem_1034}(i) we have $\Phi \left(t_0 + 2 \delta^{1/4} \right) \geq \delta / C_1$,
hence by (\ref{eq_627}),
$$ \PP
 \left( \eps \Gamma + X \geq t_0 + 2 \delta^{1/4} \right) \geq
\Phi \left(t_0 + 2 \delta^{1/4} \right) - A \eps^2 \geq \frac{\delta}{C_1} - A \eps^2
\geq \frac{\delta}{2 C_1}. $$
Consequently, from (\ref{eq_1116}),
$$ \PP \left(X \geq t_0 + \delta^{1/4} \right) \geq \PP
 \left( \eps \Gamma + X \geq t_0 + 2 \delta^{1/4} \right) - \PP \left(
\eps \Gamma \geq \delta^{1/4} \right) \geq \delta / (4 C_1). $$
A similar argument, using Lemma \ref{lem_1034}(ii) in place of Lemma \ref{lem_1034}(i), shows that  $$ \PP \left(X \leq t_0 - \delta^{1/4} \right)
\geq
\PP
 \left( \eps \Gamma + X \leq t_0 - 2 \delta^{1/4} \right) - \PP \left(
|\eps \Gamma| \geq \delta^{1/4} \right) \geq \delta / (4 C_1). $$ We conclude that
for any $t \in [t_0 - \delta^{1/4}, t_0 + \delta^{1/4}]$,
\begin{equation}
 \min \left \{ \PP \left(X \geq t \right), \PP \left( X \leq t \right) \right \}
\geq \frac{\delta}{4 C_1}. \label{eq_1030}
\end{equation}

\medskip
{\bf Step 3:} The density $f$ is differentiable  and positive everywhere. Fix $x_0 \in \RR$. Since $\log f$ is concave, then
$$
 f(x) \leq f(x_0) \exp \left( \frac{f^{\prime}(x_0)}{f(x_0)}  (x - x_0) \right)
\ \ \ \ \ \ \ \ \forall x \in \RR.
$$
Consequently, when $f^{\prime}(x_0) \neq 0$,
\begin{eqnarray*}
\lefteqn{ \min \left \{ \int_{x_0}^{\infty} f(x) dx,
\int_{-\infty}^{x_0} f(x) dx \right \} } \\ & \leq & \int_{x_0}^{\infty} f(x_0) \exp \left(
- \frac{|f^{\prime}(x_0) (x - x_0)|}{f(x_0)}  \right) dx = \frac{f(x_0)^2}{|f^{\prime}(x_0)|}.
\end{eqnarray*}
We conclude from (\ref{eq_1030}) that for any $t \in [t_0 - \delta^{1/4}, t_0 + \delta^{1/4}]$,
\begin{equation} |f^{\prime}(t)| \leq  f^2(t) \left[  \min \left \{ \PP \left(X \geq t \right), \PP \left( X \leq t \right) \right \} \right]^{-1} \leq 4 C_1 \delta^{-1} f^2(t).
\label{eq_1045}
\end{equation}
Equivalently, $|(1 / f)^{\prime}| \leq 4 C_1 \delta^{-1}$ in the
interval $[t_0 - \delta^{1/4}, t_0 + \delta^{1/4}]$. Hence,
$$ \left| \frac{1}{f(t)} - \frac{1}{f(t_0)} \right| \leq 4 C_1 \delta^{-1}
\cdot \frac{c_2}{4 C_1} \delta^{1/4} = c_2 \delta^{-3/4}
\ \ \ \ \ \ \text{when} \ \ \ |t  - t_0| \leq \frac{c_2}{4 C_1} \delta^{1/4}, $$
for $c_2 > 0$ being the universal constant from Lemma \ref{lem_1034}.
Recall from (\ref{eq_114}) that $f(t_0) = \vphi(t_0)$.
Lemma \ref{lem_1034}(iii) thus  implies that
$$ f^2(t) \leq C_1 \delta \ \ \ \ \ \ \ \ \ \text{for} \ \ \ t \in [t_0 - c \delta^{1/4}, t_0 + c \delta^{1/4}], $$
with $c = c_2 / 4 C_1$. Returning to (\ref{eq_1045}), we finally deduce the bound
$$
 |f^{\prime}(t)| \leq \tilde{C}
\ \ \ \ \ \ \ \ \ \text{for} \ \ \ t \in [t_0 - \hat{c} \delta^{1/4}, t_0 + \hat{c} \delta^{1/4}].
$$
Through Taylor's theorem, the latter bound entails that
\begin{equation}
\PP(X \geq t_0 + s) = \PP(X \geq t_0) - f(t_0)  s
+ O \left( s^2 \right)
 \ \ \ \text{for any} \ |s| \leq \hat{c} \delta^{1/4}.
\label{eq_1108}
\end{equation}

\medskip
{\bf Step 4:}  Let $\eta: \RR \rightarrow [0, \infty)$ stand for the probability density of $\eps \Gamma$. The function $\eta$ is even. Recall that $\delta \geq \eps^2$. Hence,
\begin{equation}
 \int_{|s| \geq \hat{c} \delta^{1/4}} \eta(s) ds = \PP \left( |\eps \Gamma| \geq \hat{c} \delta^{1/4} \right)
\leq \frac{\eps^4 \EE \Gamma^4}{ \hat{c}^4 \delta} \leq C \eps^2, \label{eq_540} \end{equation}
 where $\hat{c} > 0$ is the constant from (\ref{eq_1108}).
The crucial observation is that  $s \mapsto f(t_0) s \eta(s)$
is an odd function, hence its integral on a symmetric interval about the origin vanishes.
By (\ref{eq_1108}) and (\ref{eq_540}),
\begin{eqnarray*}
\lefteqn{ \left| \PP (\eps \Gamma + X \geq t_0) - \PP (X \geq t_0) \right| } \\ & = & \left| \int_{-\infty}^{\infty}
\left[ \, \PP \left( X \geq t_0 + s \right) - \PP \left(X \geq t_0 \right) \, \right] \eta(s) ds \right|
\\ & \leq & \left| \int_{-\hat{c} \delta^{1/4}}^{\hat{c} \delta^{1/4}} \left[
\,  - f(t_0) s
+ O \left( s^2 \right)  \, \right]  \eta(s) ds \right| \, + \,
2 \int_{|s| \geq \hat{c} \delta^{1/4}} \eta(s) ds
   \\ & \leq  &
 \bar{C} \int_{-\hat{c} \delta^{1/4}}^{\hat{c} \delta^{1/4}} s^2 \eta(s) ds + C \eps^2 \leq
 \bar{C} \EE  (\eps \Gamma)^2 + C \eps^2 \leq
\check{C} \eps^2,
\end{eqnarray*}
where $\hat{c} > 0$ is the constant from (\ref{eq_1108}). We apply (\ref{eq_627}) and conclude that
$$ E(t_0) = \left| \PP(X \geq t_0) - \Phi(t_0) \right| \leq \check{C} \eps^2 +
\left| \PP ( \eps \Gamma + X \geq t_0) - \Phi(t_0) \right| \leq \check{C} \eps^2 + A \eps^2. $$
Since $E(t_0) = \max_t E(t)$, the proof of the lemma is complete. \hfill $\square$

\medskip {\it Proof of Theorem \ref{thm_511}:} Let $\theta_1,\ldots,\theta_n \in \RR$ be such that $\sum_i \theta_i^2 = 1$. 
Denote $\eps = 10 \sqrt{\sum_{i=1}^n \theta_i^4}$. According
to Lemma \ref{lem_612}, the random variable
$Y = \sum_{i=1}^n \theta_i X_i$ satisfies
\begin{equation}
 \sup_{t \in \RR}
\left| \PP \left(   \eps \Gamma + Y \geq t \right)  - \Phi(t) \right| \leq C \eps^2,
\label{eq_1133} \end{equation}
with some universal constant $C \geq 1$. The random variable $Y$ has an even, log-concave density by Pr\'ekopa-Leindler.
We may thus apply Lemma \ref{lem_618}, and conclude from (\ref{eq_1133}) that
$$ \sup_{\alpha \leq \beta} \left| \PP \left(  \alpha \leq Y \leq \beta \right)
 - \left[ \Phi(\alpha) - \Phi(\beta) \right] \right| \leq 2 \sup_{t \in \RR}
\left| \PP \left(  Y \geq t \right)  - \Phi(t) \right| \leq C^{\prime} \eps^2. $$
The theorem is thus proven. \hfill $\square$

\bigskip \noindent {\large \bf Appendix: Proof of Theorem
\ref{thm_258}} \bigskip

With C\'edric Villani's permission, we reproduce below
the proof of Theorem \ref{thm_258} from his book
 \cite[Section 7.6]{villani} with a few minor changes.

\medskip
\emph{Proof of Theorem \ref{thm_258}:} We need to prove that for any $C^{\infty}$-smooth
function $\vphi: \RR^n \rightarrow \RR$,
\begin{equation}
 \int_{\RR^n} h \vphi d \mu \leq \sqrt{ \int_{\RR^n} |\nabla \vphi|^2  d \mu } \cdot
 \liminf_{\eps \rightarrow 0^+} \frac{W_2(\mu, \mu_{\eps})}{\eps}.
\label{eq_303}
\end{equation}
Since $\mu$ is compactly-supported, it is enough to restrict
attention to compactly-supported functions $\vphi$. Fix such a
test function $\vphi$. Then the second derivatives of $\vphi$ are
bounded on $\RR^n$. By Taylor's theorem, there exists a constant $R = R(\vphi)$
with
\begin{equation}
\vphi(y) - \vphi(x) \leq | \nabla \vphi(x) | \cdot |x - y| + R
|x-y|^2 \ \ \ \ \ \ \forall x,y \in \RR^n. \label{eq_316}
\end{equation}
We may assume that $\sup |h| > 0$ (otherwise, the theorem holds
trivially), and let $\eps
> 0$ be smaller than $1 / \sup |h|$. Then
$\mu_{\eps}$ is a non-negative measure on $\RR^n$. Let $\gamma$ be
any coupling of $\mu$ and $\mu_{\eps}$. We see that
$$ \int_{\RR^n} h \vphi d \mu = \frac{1}{\eps} \int_{\RR^n}
 \vphi d \left[ \mu_{\eps} - \mu \right]
= \frac{1}{\eps} \int_{\RR^n \times \RR^n} \left[ \vphi(y) -
\vphi(x) \right] d \gamma(x,y). $$ Write $W_2^{\gamma}(\mu,
\mu_{\eps}) = \sqrt{ \int_{\RR^n \times \RR^n} |x - y|^2 d
\gamma(x,y) }$. According to (\ref{eq_316}) and to the
Cauchy-Schwartz inequality,
\begin{eqnarray*}
\int_{\RR^n}  h \vphi d \mu & \leq & \frac{1}{\eps} \int_{\RR^n \times \RR^n} |\nabla \vphi(x)| \cdot |x - y| d \gamma(x,y) + \frac{R}{\eps} \int_{\RR^n \times \RR^n} |x - y|^2 d \gamma(x,y)  \\
 & \leq & \frac{1}{\eps} \sqrt{ \int_{\RR^n} |\nabla \vphi(x)|^2 d \mu(x) } \cdot W_2^{\gamma}(\mu, \mu_{\eps})
+ \frac{R}{\eps} W_2^{\gamma}(\mu, \mu_{\eps})^2.
\end{eqnarray*}
By taking the infimum over all couplings $\gamma$ of $\mu$ and
$\mu_{\eps}$, we obtain
\begin{equation}
 \int_{\RR^n}  h \vphi d \mu
 \leq \sqrt{ \int_{\RR^n} |\nabla \vphi|^2 d \mu } \cdot \frac{W_2(\mu, \mu_{\eps})}{\eps}
+ R \frac{W_2(\mu, \mu_{\eps})^2}{\eps} , \label{eq_327}
\end{equation}
with $R$ depending only on $\vphi$.
We may assume that $\liminf_{\eps \rightarrow 0^+} W_2(\mu,
\mu_{\eps}) / \eps < \infty$; otherwise, there is nothing to
prove. Consequently,
$$ \liminf_{\eps \rightarrow 0^+} \frac{W_2(\mu, \mu_{\eps})^2}{\eps} =
\liminf_{\eps \rightarrow 0^+} \eps \left( \frac{W_2(\mu,
\mu_{\eps})}{\eps} \right)^2 = 0. $$ Hence by letting $\eps$ tend
to zero in (\ref{eq_327}), we deduce (\ref{eq_303}). The proof is
complete. \hfill $\square$


\end{document}